\documentclass{amsart}
\usepackage{hyperref}
\usepackage[table,dvipsnames]{xcolor}
\usepackage{color,amssymb, comment, mathrsfs, tikz, upgreek, contour, soul, cleveref, mathtools, enumerate, pdfpages, cleveref, adjustbox, tikz-cd}
\usetikzlibrary{cd,fit,shapes.geometric,matrix,arrows}

\theoremstyle{definition}

\theoremstyle{remark}

\newcommand{\ee}{\mathsf{e}}
\newcommand{\ff}{\mathsf{f}}
\newcommand{\ZZ}{\mathbb{Z}}
\newcommand{\QQ}{\mathbb{Q}}
\renewcommand{\gg}{\mathsf{g}}

\DeclareMathOperator{\Tor}{\mathrm{Tor}}
\DeclareMathOperator{\Hom}{\mathrm{Hom}}
\DeclareMathOperator{\rank}{\mathrm{rank}}
\renewcommand{\H}{\mathbf{H}}
\DeclareMathOperator{\B}{\mathbf{B}}
\DeclareMathOperator{\G}{\mathbf{G}}
\DeclareMathOperator{\T}{\mathbf{T}}

\DeclareMathOperator{\fieldChar}{\mathsf{fieldChar}}
\DeclareMathOperator{\logging}{\mathsf{logging}}
\DeclareMathOperator{\numTerms}{\mathsf{numTerms}}
\DeclareMathOperator{\mn}{\mathsf{mn}}
\DeclareMathOperator{\lowDeg}{\mathsf{lowDeg}}
\DeclareMathOperator{\highDeg}{\mathsf{highDeg}}
\DeclareMathOperator{\degSeq}{\mathsf{degSeq}}
\DeclareMathOperator{\checkIn}{\mathsf{checkIn}}
\DeclareMathOperator{\strictTerms}{\mathsf{strictTerms}}
\DeclareMathOperator{\maxTries}{\mathsf{maxTries}}
\DeclareMathOperator{\useN}{\mathsf{useN}}
\DeclareMathOperator{\maxM}{\mathsf{maxM}}
\DeclareMathOperator{\maxN}{\mathsf{maxN}}

\numberwithin{equation}{section}

\allowdisplaybreaks

\begin{document}

\title[Sampling Algebra Structures on Minimal Free
Resolutions]{Sampling Algebra Structures\\ on Minimal Free Resolutions
  of length $\mathbf{3}$}

\author[L.W.~Christensen]{Lars Winther Christensen} \thanks{L.W.C. was
  partly supported by Simons Foundation collaboration grants 428308
  and 962956} \address{Department of Mathematics, Texas Tech
  University, Lubbock, TX 79409,
  U.S.A.}\email{lars.w.christensen@ttu.edu} \author[O.~Gotchey]{Orin
  Gotchey} \address{Department of Mathematics, Texas Tech University,
  Lubbock, TX 79409,
  U.S.A.}\email{orin.gotchey@ttu.edu}\curraddr{orin.gotchey@skyeanalytics.com}
\author[A.~Hardesty]{Alexis Hardesty} \address{Division of
  Mathematics, Texas Woman's University, Denton, TX 76204, U.S.A.}
\email{ahardesty1@twu.edu}



\date{25 September 2024}

\begin{abstract}
  Ideals in the ring of power series in three variables over a field
  can be classified based on algebra structures on their minimal free
  resolutions. The classification is incomplete in the sense that it
  remains open which algebra structures actually occur; this
  \emph{realizability question} was formally raised by Avramov in
  2012. We discuss the outcomes of an experiment performed to shed
  light on Avramov's question: Using the computer algebra system
  \emph{Macaulay2}, we classify a billion randomly generated ideals
  and build a database with examples of ideals of all classes realized
  in the experiment. Based on the outcomes, we discuss the status of
  conjectures that relate to the realizability question.
\end{abstract}

\maketitle

\section*{Introduction}

\noindent
Let $R$ be a local ring with residue field $\Bbbk$ and $I \subset R$ a
perfect ideal of grade $3$. By a result of Buchsbaum and Eisenbud
\cite[Proposition 1.1]{BE77}, the minimal free resolution $F_\bullet$
of $R/I$ over $R$ has a differential graded (DG) algebra
structure. This induces a graded $\Bbbk$-algebra structure on
$\mathrm{H}(F_\bullet\otimes_R\Bbbk) = \Tor_\bullet^R(R/I,\Bbbk)$, and
while the DG algeba structure on $F_\bullet$ is not unique, the
induced $\Bbbk$-algebra structure on $\Tor_\bullet^R(R/I,\Bbbk)$ is
unique. Results of Weyman \cite[Theorem 4.1]{Weyman89} and of Avramov,
Kustin, and Miller \cite[Theorem 2.1]{AKM88} show that this structure
supports a classification scheme for grade $3$ perfect ideals in
$R$. The original application of the classification scheme was to
answer a question in local algebra---on the rationality of Poincar\'e
series---and for this it was sufficient to establish the possible
structures without considering their realizability. Later, however,
Avramov returned to the classification to resolve another question in
local algebra---on growth patterns in minimal injective
resolutions---and found it necessary to rule out the realizability of
certain structures. Ideally, a classification should say exactly which
structures occur, and Avramov formally stated that question in
\cite[Question 3.8]{Avramov12}.

Since \cite{Avramov12}, various authors have addressed what has become
known as the \emph{realizability question} in essentially two
different ways: Some have ruled out the realizability of certain
classes while others have provided constructions of ideals in certain,
other, classes. In this paper the approach is experimental: Within
certain bounds we generated random grade $3$ perfect ideals and
classified them. Then we compare our observations to existing bounds,
both established and conjectured.  In addition to the results and
analysis in this paper, we provide a \textit{GitHub} repository with
the simplest examples of ideals from each class that was observed in
the experiment. It turns out that many classes can be realized by
binomial ideals and some even by monomial ideals.

\begin{equation*}
  \ast \ \ast \ \ast    
\end{equation*}

Set $A_\bullet=\Tor_\bullet^R(R/I,\Bbbk)$. The size of the algebra
$A_\bullet$ is determined by two parameters, $m$ and $n$, where $m$ is
the minimal number of generators of the ideal $I$ and $n$ is known as
the \emph{type} of the ring $R/I$. It is proved in \cite[Theorem
2.1]{AKM88} that there exists bases
\[
  \{\ee_i\}_{i=1,\ldots,m},\quad\{\ff_i\}_{i=1,\ldots,m+n-1},\quad\textnormal{and}\quad\{\gg_i\}_{i=1,\ldots,n}
\]
for $A_1,A_2$ and $A_3$, respectively, such that the multiplication on
$A_\bullet$ is one of the following:

\smallskip
\begin{center}
  \renewcommand{\arraystretch}{1.1}
  \begin{tabular}{rlrl}
    $\mathbf{C}(3)$ & $\ee_1 \ee_2 = \ff_3$, $\ee_2 \ee_3 = \ff_1$, $\ee_3 \ee_1 = \ff_2$ & $\ee_i \ff_i = \gg_1$ & for $1\leq i\leq 3$\\
    $\T$ & $\ee_1 \ee_2 = \ff_3$, $\ee_2 \ee_3 = \ff_1$, $\ee_3 \ee_1 = \ff_2$\\
    $\B$ & $\ee_1 \ee_2 = \ff_3$ & $\ee_i \ff_i = \gg_1$ & for $1\leq i \leq 2$\\
    $\G(r)$ &  & $\ee_i \ff_i = \gg_1$ & for $1\leq i\leq r$\\
    $\H(p,q)$ & $\ee_i \ee_{p+1} = \ff_i$ for $1\leq i\leq p$ & $\ee_{p+1}\ff_{p+j} = \gg_j$ & for $1\leq j\leq q$,
  \end{tabular}
\end{center}
\smallskip

\noindent The products not listed are either zero or can be deduced
from the ones listed by graded-commutativity. Set
\begin{equation*}
  m=\rank A_1, \quad n=\rank A_3, \quad p=\rank A_1A_1, \quad q=\rank
  A_1A_2 \textnormal{ and } r=\rank\delta
\end{equation*}
where $\delta$ is the canonical homomorphism
$A_2 \rightarrow \Hom_\Bbbk(A_1,A_3)$ that maps an element $a$ to the
multiplication map $A_1 \xrightarrow{a\cdot} A_3$.

The realizability question for grade $3$ perfect ideals can now be
phrased as follows: Given $m$ and $n$, which of the classes above can
be realized? For small values of $m$ and $n$ the answer is known:
Ideals of class $\mathbf{C}(3)$, i.e.\ grade $3$ complete intersection
ideals, only exist for $(m,n) = (3,1)$. For $(m,n)$ with $m \le 4$ or
$n=1$ there is at most one permissible class and every permissible
class can be realized; see Avramov \cite[1.4.2, 3.4.1(a), 3.4.2,
3.9.1]{Avramov12}. For this reason, we only recorded ideals with
$m \ge 5$ and $n \ge 2$. For practical reasons, we also only recorded
ideals with $m \le 12$ and $n \le 10$.  The results of our experiment
indicate that there are no further restrictions on the realizability
of classes $\B$, $\H$, and $\T$ beyond those proved in \cite[Theorem
3.1]{Avramov12} and \cite[Theorem 1.1]{CVW20Linkage}, but they do not
rule out that there may be restrictions on the existence of class $\G$
ideals beyond what is known or conjectured in the literature.

\section{Materials and Methods}\label{sec1}

\noindent
To shed light on the realizability question, we developed an algorithm
that randomly generates homogeneous grade $3$ perfect ideals in the
trivariate polynomial ring $\Bbbk[x,y,z]$ and classifies the
corresponding ideal in the local ring
$\Bbbk[\mspace{-2.5mu}[x,y,z]\mspace{-2mu}]$ using the
\emph{Macaulay2} package \emph{TorAlgebra}
\cite{M2-LWCOVl,LWCOVl14}. To be precise, given an ideal
$I = (f_1,\ldots,f_m)$ in the polynomial algebra, the quotient
$R = \Bbbk[x,y,z]/I$ is an Artinian local ring isomorphic to
$\Bbbk[\mspace{-2.5mu}[x,y,z]\mspace{-2mu}]/I'$, where $I'$ denotes
the codepth $3$ perfect ideal of the power series algebra
$\Bbbk[\mspace{-2.5mu}[x,y,z]\mspace{-2mu}]$ generated by
$f_1,\ldots,f_m$. The multiplicative structure on the free resolution
of the quotient ring $R$ is encoded in homological invariants of that
ring. As the classification algorithm from \cite{LWCOVl14} implemented
in \emph{TorAlgebra} classifies the local ring $R$ based on these
intrinsic properties, it is irrelevant that it is obtained as a
quotient of the (nonlocal) polynomial algebra rather than the local
ring $\Bbbk[\mspace{-2.5mu}[x,y,z]\mspace{-2mu}]$. The algorithm
maintains a running tally of the classified ideals, along with a
running list of the ``shortest" minimal generating set observed for
ideals of each class. This section describes the algorithm in detail,
highlighting the user-readable output on screen and in the recorded
data structure.

The entry point into the algorithm is a function called \textsf{main},
which takes one positional parameter and eleven optional
parameters. The positional parameter is a positive integer which
controls the number of attempts to classify ideals. The optional
parameters are detailed below.

\smallskip
\begin{center}
  \renewcommand{\arraystretch}{1.1}
  \begin{tabular}{rllll}
    Parameter && Data Type && Default Value\\
    \hline
    $\fieldChar$ && $0$ or a prime integer && 3\\
    $\checkIn$ && nonnegative integer && 0\\
    $\degSeq$ && $(0)$ or a sequence of positive integers && $(0)$\\
    $\lowDeg$ && positive integer && 2 \\
    $\highDeg$ && positive integer && 8\\
    $\numTerms$ && nonnegative integer && 0\\
    $\mn$ && positive integer && 5\\
    $\useN$ && boolean && false\\
    $\maxTries$ && positive integer && 10\\
    $\strictTerms$ && boolean && false\\
    $\maxM$ && positive integer && 12 \\
    $\maxN$ && positive integer && 10\\
    $\logging$ && boolean && false
  \end{tabular}
\end{center}
\smallskip

\noindent When the main function is called, it executes the following
steps in order:

\subsection{Set the Polynomial Ring} \label{step1} The experiment
takes place within an ambient polynomial ring, $R=\Bbbk[x,y,z]$, where
$\Bbbk$ is a field.  If $\fieldChar$ is a prime $p$, then $\Bbbk$ is
set to $\mathbb{Z}/p\ZZ$. If $\fieldChar=0$, then $\Bbbk$ is set to
$\QQ$. Otherwise, the function prints the statement
\begin{align*}
  \text{Error: bad field.}
\end{align*}

\subsection{Load Pre-existing Data} \label{step2} The function
searches the current working directory for a \texttt{data} folder from
a previous execution. If it finds one, then the function loads data
from that folder into memory. Otherwise, the function creates an empty
directory, \texttt{data}.

\subsection{Print Start Message} \label{step3} The following message
is printed to the terminal:
\begin{align*}
  &\text{Main Routine started at } \textsf{current time}  \text{ with options:}\\
  &\text{new OptionTable from } \{ 
    \text{maxTries $=> \maxTries$, degSeq $=> \degSeq$,}\\
  &\text{strictTerms $=> \strictTerms$, logging $=> \logging$, mn $=> \mn$,}\\ 
  &\text{numTerms $=> \numTerms$, highDeg $=> \highDeg$, useN $=> \useN$,}\\
  &\text{maxM $=> \maxM$, maxN $=> \maxN$, checkIn $=> \checkIn$,}\\
  &\text{fieldChar $=> \fieldChar$, lowDeg $=> \lowDeg$}\}
\end{align*}
If $\logging$ is true, then the function appends the printed statement
to the end of the \texttt{log.txt} file.

\medskip
\begin{center}
  \textbf{Steps \ref{step4} and \ref{step5} form the loop of the
    algorithm.}
    
  The positional parameter controls the number of repetitions.
\end{center}

\subsection{Generate an Ideal} \label{step4} The function checks if
the repetition counter $i$ is a multiple of $\checkIn$. If so, then
the function prints the statement
\begin{align*}
  \text{Checking in every $\checkIn$ ideals... done $i$ so far.}    
\end{align*}
If $\degSeq$ is the sequence $(0)$, then the function changes it to a
sequence of $\mn$ randomly generated integers within the interval
$[\lowDeg,\,\highDeg]$. Next, a new sequence $S$ of length $\mn$ is
created, consisting of random homogeneous polynomials of degree $d$
for each element $d$ of $\degSeq$. If $\numTerms>0$, then these
polynomials are constructed so that each has $\numTerms$ terms
(e.g. $\numTerms=1$ causes the algorithm to produce monomials).  Else,
i.e.\ if $\numTerms=0$, then the number of terms in each polynomial is
random. These polynomials are then used to generate a homogeneous
ideal as follows:

\begin{enumerate}[i.]

\item If $\useN$ is false, then the function creates an ideal $I$
  which is generated (possibly non-minimally) by the elements of
  $S$. The function checks if $I$ is minimally generated by the $\mn$
  elements of $S$, trying up to $10$ times to do so. Each time, the
  function generates an altered ideal by taking the previous
  generating set and adding a random homogeneous form of degree
  randomly selected from $\degSeq$. If after $10$ tries the function
  has still not constructed an ideal minimally generated by $\mn$
  elements, it considers the attempt a failure and restarts Step
  \ref{step4}. Otherwise:
  \begin{itemize}
  \item If $\text{codim}\,I=3$, then the ideal has been successfully
    generated, and the function moves on to Step \ref{step5}.
  \item If $\text{codim}\,I<3$ and $\numTerms\ne 1$, then the function
    adds a pure power of a variable to one or more of the minimal
    generators to construct a new ideal, keeping the polynomial
    generators homogeneous and $\mn$ as the minimal number of
    generators. If needed, this is repeated for each variable, always
    starting with the minimal set of generators. This process stops
    when an ideal of codimension 3 and $\mn$ generators is achieved,
    or when the set of variables is exhausted. In the latter case, the
    $0$ ideal is returned.
  \item If $\text{codim}\,I<3$ and $\numTerms=1$, then the function is
    intended to return a monomial ideal. In this case, the above
    strategy to fix $\text{codim}\,I$ is not appropriate. Therefore,
    the function considers this attempt a failure and restarts Step
    \ref{step4}.
  \end{itemize}
   
\item If $\useN$ is true, then the function attempts to create an
  ideal with quotient of type $\mn$. It applies the Macaulay2 function
  $\textsf{fromDual}$ to the elements of the sequence $S$, which
  returns a set of generators of a homogeneous ideal that defines a
  quotient ring of type (at most) $\mn$; see for example Meyer and
  Smith \cite[Chapter II.2]{pms}. If the type is $\mn$, then the
  algorithm has succeeded in generating an ideal with the required
  properties. If not, then the function considers this attempt a
  failure and restarts Step \ref{step4}.
\end{enumerate}

\noindent If the function considers an attempt a failure, then a
variable called \textit{numTries} (initially set to $0$) is checked
against the optional variable $\maxTries$. If
$\textit{numTries}<\maxTries$, then the former is incremented before
the function restarts Step \ref{step4}. On the other hand, if
$\textit{numTries}=\maxTries$, then the function returns the $0$
ideal.  If the function succeeds in generating an ideal with the
required properties, then \textit{numTries} is reset to $0$.

If $\strictTerms$ is true, then the function checks if the minimal
generators of the ideal have the exact number of terms as given by
$\numTerms$. In addition, the function verifies that the ideal has
codimension 3, that it is homogeneous, that its minimal number of
generators does not exceed $\maxM$, that the type of its quotient does
not exceed $\maxN$, and that all minimal generators are of degree at
least $2$. If any of these checks fail, then the ideal is not
classified.

\subsection{Classify the Ideal and File Classification
  Data} \label{step5} The function classifies the ideal using the
TorAlgData command in the \textit{TorAlgebra}\cite{M2-LWCOVl} package,
resulting in a tuple $(m,n,\textnormal{Class},p,q,r)$. If they do not
already exist, then the function creates the following files with the
\texttt{data} folder: \texttt{classDat.txt} and \texttt{class.txt},
the former being \emph{Macaulay2} readable and the latter intended to
be human readable. If $\numTerms=0$, then the function edits files in
the subfolder $\texttt{data/0}$---creating it if necessary. Otherwise,
the function computes the maximum number of terms of a minimal
generator of the ideal and edits files in a subfolder---creating it if
necessary---called \texttt{1}, \texttt{2}, \texttt{3}, or \texttt{4}
within the \texttt{data} folder corresponding to monomial, binomial,
trinomial, or generators with four or more terms.

\begin{enumerate}[i.]
\item If the class has not been seen before, then the function creates
  a \texttt{.txt} file named \texttt{m-n-Class-p-q-r} with a
  \emph{Macaulay2} readable matrix containing the minimal generators
  of the ideal. For example, the file \texttt{5-2-B-1-1-2.txt} in
  \texttt{data/2} could have the contents:

  \smallskip
  \begin{center}
    \verb|matrix{{y*z,x*z,y^2+z^2,x*y+z^2,x^2+z^2}}|
  \end{center}
  \smallskip
  
  \noindent In the \texttt{data} folder, the function adds the class
  to the \texttt{classDat.txt} file, recording the tuple
  $(m,n,\textnormal{Class},p,q,r)$, the \emph{Macaulay2} readable
  matrix containing the minimal generators of the ideal, and a count
  corresponding to the number of times the class has been observed. In
  the running example, the corresponding entry in the file
  \texttt{classDat.txt} would be

  \smallskip
  \begin{center}
    \verb|((5,2,B,1,1,2),(matrix{{y*z,x*z,y^2+z^2,x*y+z^2,x^2+z^2}},1))|
  \end{center}
  \smallskip
  
  \noindent Finally, in the \texttt{data} folder, the function adds a
  row to the \texttt{class.txt} file with the same information,
  replacing the \emph{Macaulay2} readable matrix containing the
  minimal generators of the ideal with a human readable list. Note
  that the \texttt{class.txt} file is ordered numerically according to
  the value of $m$, then ordered numerically according to the value of
  $n$, then ordered alphabetically according to Class, then ordered
  numerically according to the value of $p$, then $q$, then~$r$. In
  the running example, the first entry in the \texttt{class.txt} file
  would be
  
  \smallskip
  \begin{center}
    \verb#| 5 2 B 1 1 2 1 | yz xz y2+z2 xy+z2 x2+z2 |#
  \end{center}
  \smallskip

\item If the class has been seen before, then the function opens the
  previously created \texttt{m-n-Class-p-q-r.txt} file and compares
  the length of the minimal generators of the current ideal to the
  length of the previously recorded minimal generators. Here ``length"
  simply refers to the length of the text string. If the length of the
  minimal generators of the current ideal is shorter than the length
  of the previously recorded minimal generators, then a
  \emph{Macaulay2} readable matrix containing the minimal generators
  of the current ideal is appended to the end of the file. For
  example, an updated \texttt{5-2-B-1-1-2.txt} file in \texttt{data/2}
  could have the contents:

  \smallskip
  \begin{center}
    \begin{tabular}{l}
      \verb|matrix{{y*z,x*z,y^2+z^2,x*y+z^2,x^2+z^2}}|\\
      \verb|matrix{{z^2,y*z,x*z,x*y,x^2-y^2}}|
    \end{tabular}
  \end{center}
  \smallskip
  
  \noindent Additionally, the function increases the count of the
  class in the \texttt{classDat.txt} and \texttt{class.txt} files by l
  and replaces the previously recorded minimal generators of the ideal
  with the minimal generators of the current ideal, provided their
  length is shorter. In the running example, the updated entries in
  \texttt{classDat.txt} and \texttt{class.txt} would be

  \smallskip
  \begin{center}
    \verb|((5,2,B,1,1,2),(matrix{{z^2,y*z,x*z,x*y,x^2-y^2}},2))|
  \end{center}
  \smallskip

  \noindent and

  \smallskip
  \begin{center}
    \verb#| 5 2 B 1 1 2 2 | z2 yz xz xy x2-y2 |#
  \end{center}
\end{enumerate}

\medskip
\noindent The function repeats Steps \ref{step4} and \ref{step5}
according to the positional parameter entered by the user when the
main function was called.

\subsection{Print Summary} \label{step6} The function prints the
following information:
\begin{align*}
  &\text{Main Routine finished:} \\
  &\text{at \textsf{current time}}\\
  &\text{ran for \texttt{\#} seconds,}\\
  &\text{classified \texttt{\#} ideals,}\\
  &\text{generated \texttt{\#} distinct classes,}\\
  &\text{discovered \texttt{\#} new classes}
\end{align*}
If new classes were discovered, then the function also prints a list
with entries
\begin{center}
  $\{(m,n,\textnormal{Class},p,q,r), \dots\}$
\end{center}
If the user set $\logging$ as true, then the printed statements above
are appended at the end of the \texttt{log.txt} file.

\section{Results}

\subsection{Realized Classes}
The total number of ideals classified in our experiment is just above
$10^9$. After some initial experimentation with the field
characteristic, we chose to work in characteristic $3$: For
computational economy the characteristic needs to be low, and compared
to characteristic $2$, ideals generated in characteristic $3$ realize
a wider variety of classes, presumably because of the existence of a
sign. It turns out that, in most cases, a set of homogeneous
polynomials with coefficients $\pm 1$ generate ideals of the same
class in $\ZZ_3[x,y,z]$ and $\QQ[x,y,z]$. Now this is only in about
85\,\% of cases; for example, the six polynomials
\begin{gather*}
  xy^2+y^2z+xz^2-z^3, \\
  x^2z^2+xyz^2-y^2z^2-yz^3-z^4, \\
  y^3z+xyz^2+y^2z^2+xz^3-yz^3, \\
  x^2yz+xz^3+z^4,\\
  x^3z+xz^3+z^4,\;\text{and}\\
  x^4+y^4+y^2z^2+xz^3+z^4
\end{gather*}
generate an ideal in $\ZZ_3[x,y,z]$ described by the tuple
$(6,4,\H,1,2,2)$ but in $\QQ[x,y,z]$ an ideal described by
$(6,6,\H,1,1,1)$.

In the experiment, the optional parameters ranged as follows:

\smallskip
\begin{center}
  \renewcommand{\arraystretch}{1.1}
  \begin{tabular}{rllll}
    Parameter && Values\\
    \hline
    $\fieldChar$ && $3$\\
    $\degSeq$ && $(2,\dots,2,2),(2,\dots,2,3),\dots,(10,\dots,10,10)$ \\
    $\lowDeg$ && 2--6 \\
    $\highDeg$ && $(\lowDeg + 0)$--$(\lowDeg + 9)$\\
    $\numTerms$ && 0--12\\
  \end{tabular}
\end{center}
\smallskip

\noindent The variation of the parameters was informed by the
outcomes: For example, when a new class was observed, a minimal set of
generators for the ideal that realized the class was recorded. These
generating sets were analyzed for patterns in their degrees and number
of terms, and the parameters were adjusted accordingly with the goal
of observing as large a variety of classes as possible. In addition,
we switched $\useN$ between true and false to observe classes of a
specific value of $n$ or $m$, respectively.

We visualize the observed classes in six tables:

\smallskip
\begin{center}
  \renewcommand{\arraystretch}{1.1}
  \begin{tabular}{rl}
    Table \ref{tableBGT} & Ideals of class $\B$, $\G$, and $\T$\\
    Table \ref{tableBGTMon} &  Monomial ideals of class $\B$, $\G$, and $\T$\\
    Table \ref{tableBGTBin} &  Binomial ideals of class $\B$, $\G$, and $\T$\\
    Table \ref{tableH} &  Ideals of class $\H$\\
    Table \ref{tableHMon} & Monomial ideals of class $\H$\\
    Table \ref{tableHBin} & Binomial ideals of class $\H$
  \end{tabular}
\end{center}
\smallskip

\noindent
For reasons of space, these tables are limited to the ranges
$5\leq m\leq 9$ and $2\leq n\leq 9$. Full tables with ranges
$5 \leq m \le 12$ and $2 \leq n \leq 10$ are available online, see
\ref{Repository}. For a fixed pair $(m,n)$, call the collection of all
classes with these $m$ and $n$ values the $(m,n)$-box.  Within each
$(m,n)$-box, the tables have either $(p,r)$-cells (Tables
\ref{tableBGT}--\ref{tableBGTBin}) or $(p,q)$-cells (Tables
\ref{tableH}--\ref{tableHBin}), according to the possible values of
$p$, $q$, and $r$, which are known to be bounded by functions of $m$
and $n$; see \cite[Theorem 1.1]{CVW20Linkage} for class $\H$ and
\cite[Theorem 3.1]{Avramov12} for class $\G$.  Dotted cells represent
classes that are known to be unrealizable; they are separated from
cells representing permissible classes with a black border. Of the
permissible classes, those that have not been observed are represented
by white boxes. Classes that have been observed in the experiment have
cells colored a shade of gray (Tables \ref{tableBGT} and \ref{tableH})
or black (Tables \ref{tableBGTMon}, \ref{tableBGTBin},
\ref{tableHMon}, and \ref{tableHBin}). Tables \ref{tableBGT} and
\ref{tableH} are colored according to the frequency with which each
class was observed in the experiments. The darker the coloring, the
more frequently the class was observed. These tables are colored using
the same scale.
\clearpage

\begin{center}
  \refstepcounter{subsubsection}
  \includegraphics[angle=90, width=4.9in]{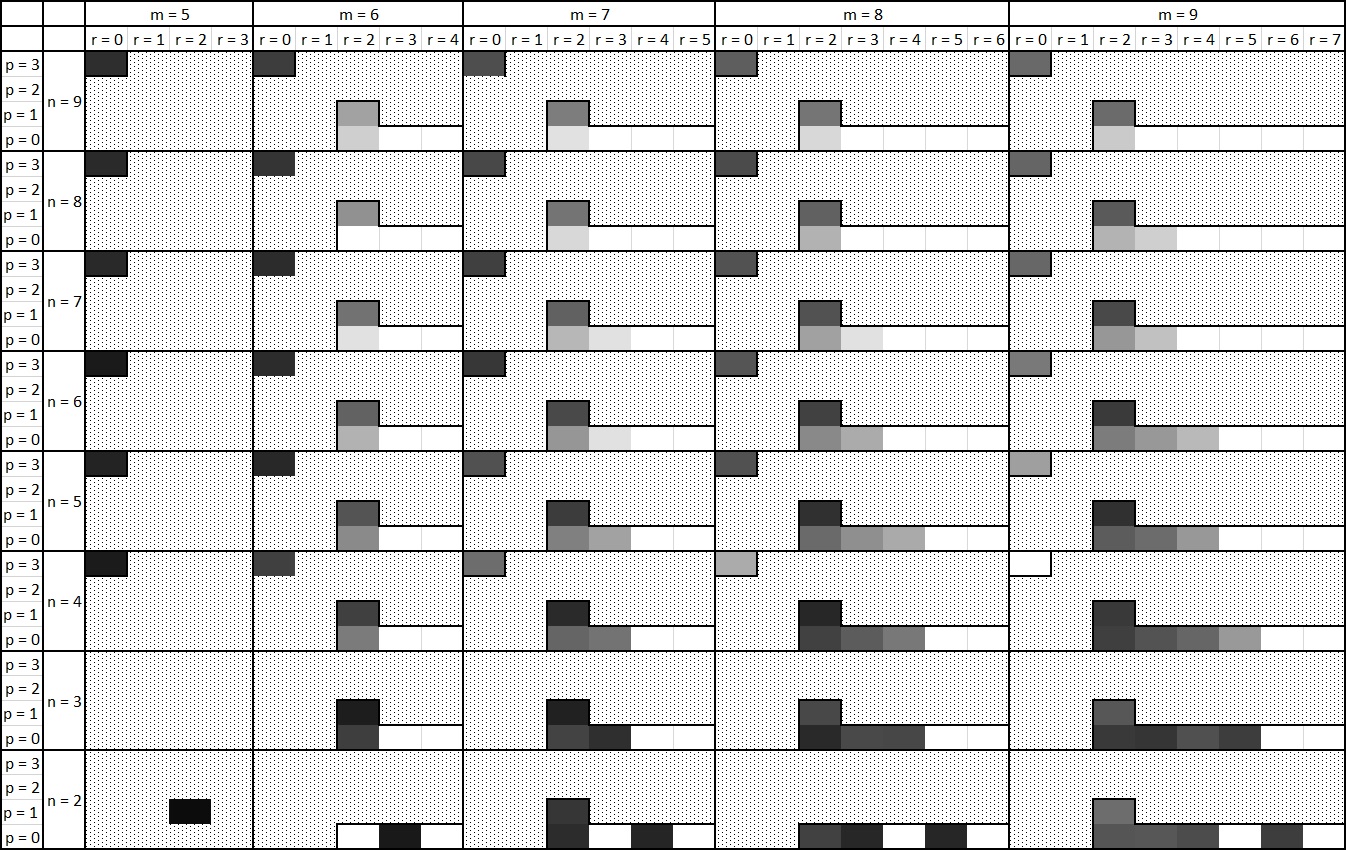}\\
  \label{tableBGT}
  Table \thesubsubsection: Observed ideals of class $\B$, $\G$, and
  $\T$.
\end{center}

\begin{center}
  \refstepcounter{subsubsection}
  \includegraphics[angle=90, width=4.9in]{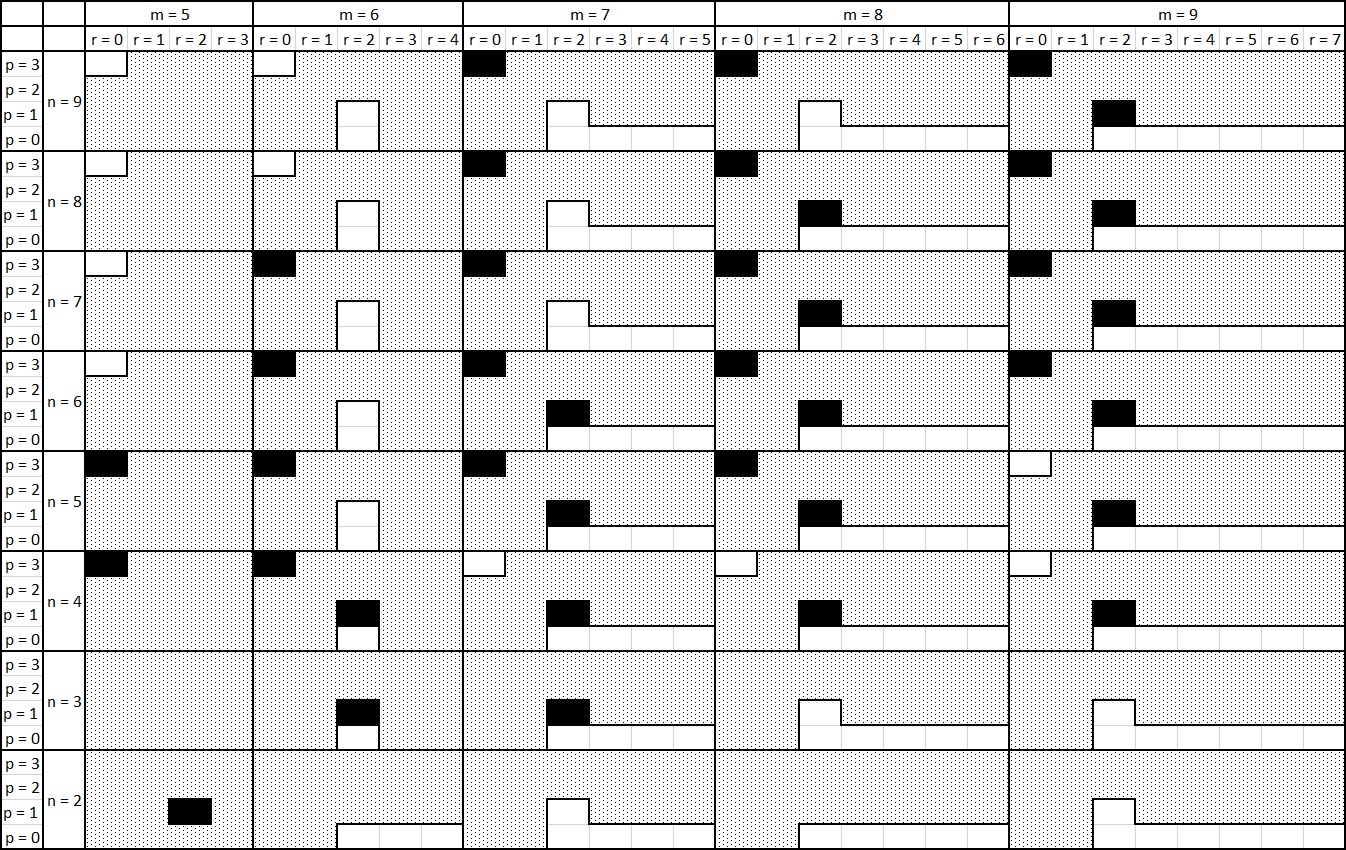}\\
  \label{tableBGTMon}
  Table \thesubsubsection: Observed monomial ideals of class $\B$,
  $\G$, and $\T$.
\end{center}

\begin{center}
  \refstepcounter{subsubsection}
  \includegraphics[angle=90, width=4.9in]{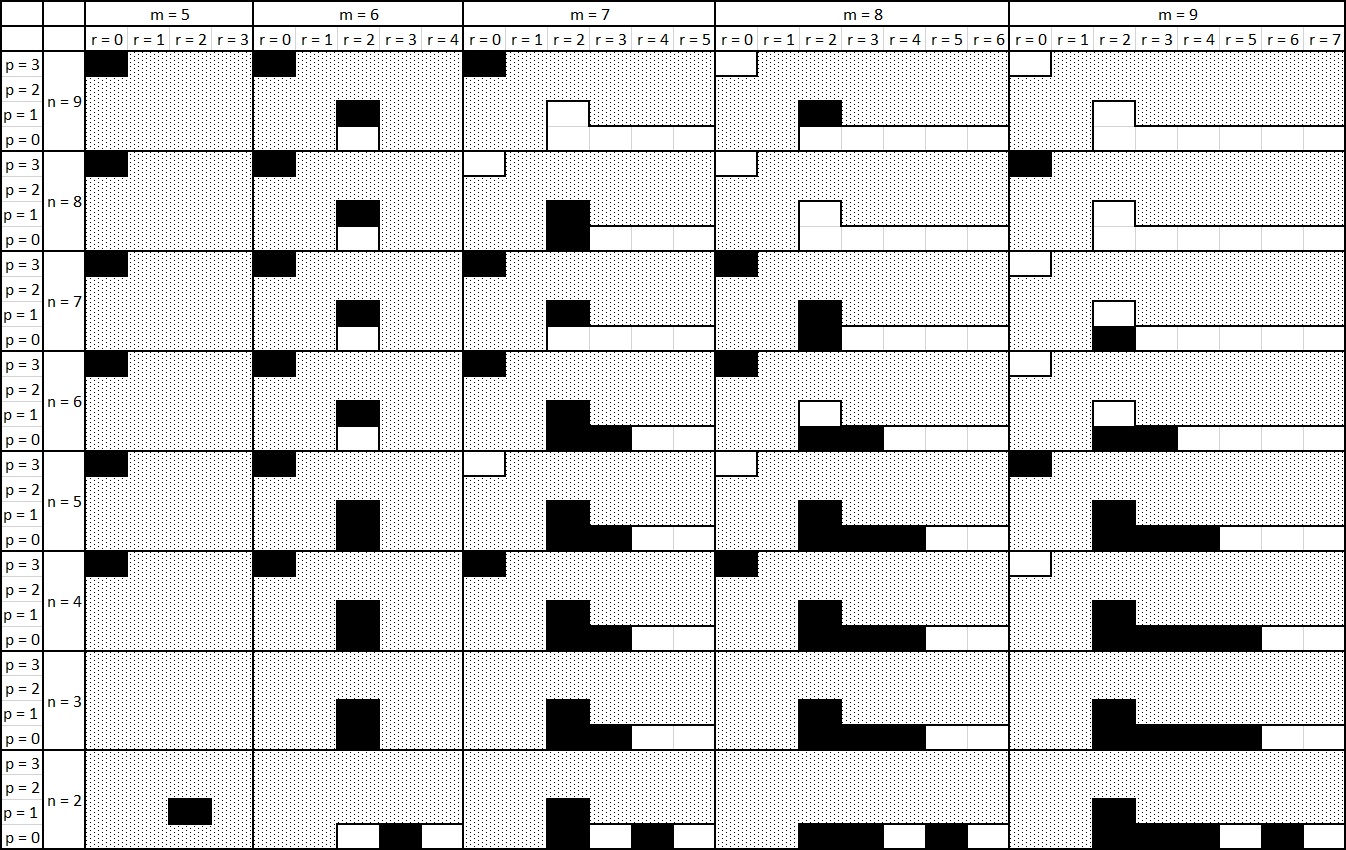}\\
  \label{tableBGTBin}
  Table \thesubsubsection: Observed binomial ideals of class $\B$,
  $\G$, and $\T$.
\end{center}

\hspace{-3.5pc}
\begin{minipage}{35pc}
  \refstepcounter{subsubsection}
  \begin{center}
    \includegraphics[angle=90, width=5.5in]{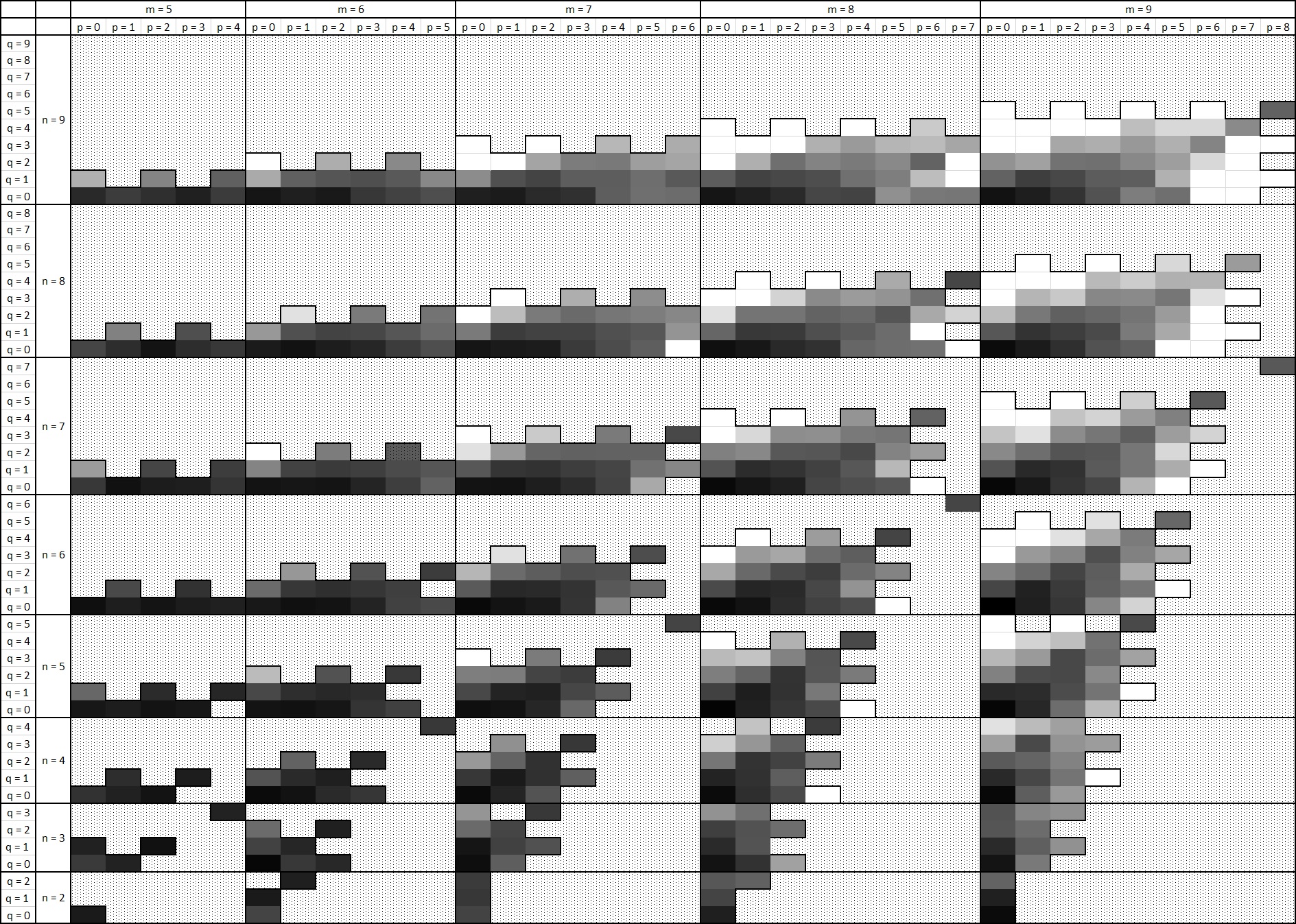}\\
    \label{tableH}
    Table \thesubsubsection: Observed ideals of class $\H$.
  \end{center}
\end{minipage}

\hspace{-3.5pc}
\begin{minipage}{35pc}
  \refstepcounter{subsubsection}
  \begin{center}
    \includegraphics[angle=90, width=5.5in]{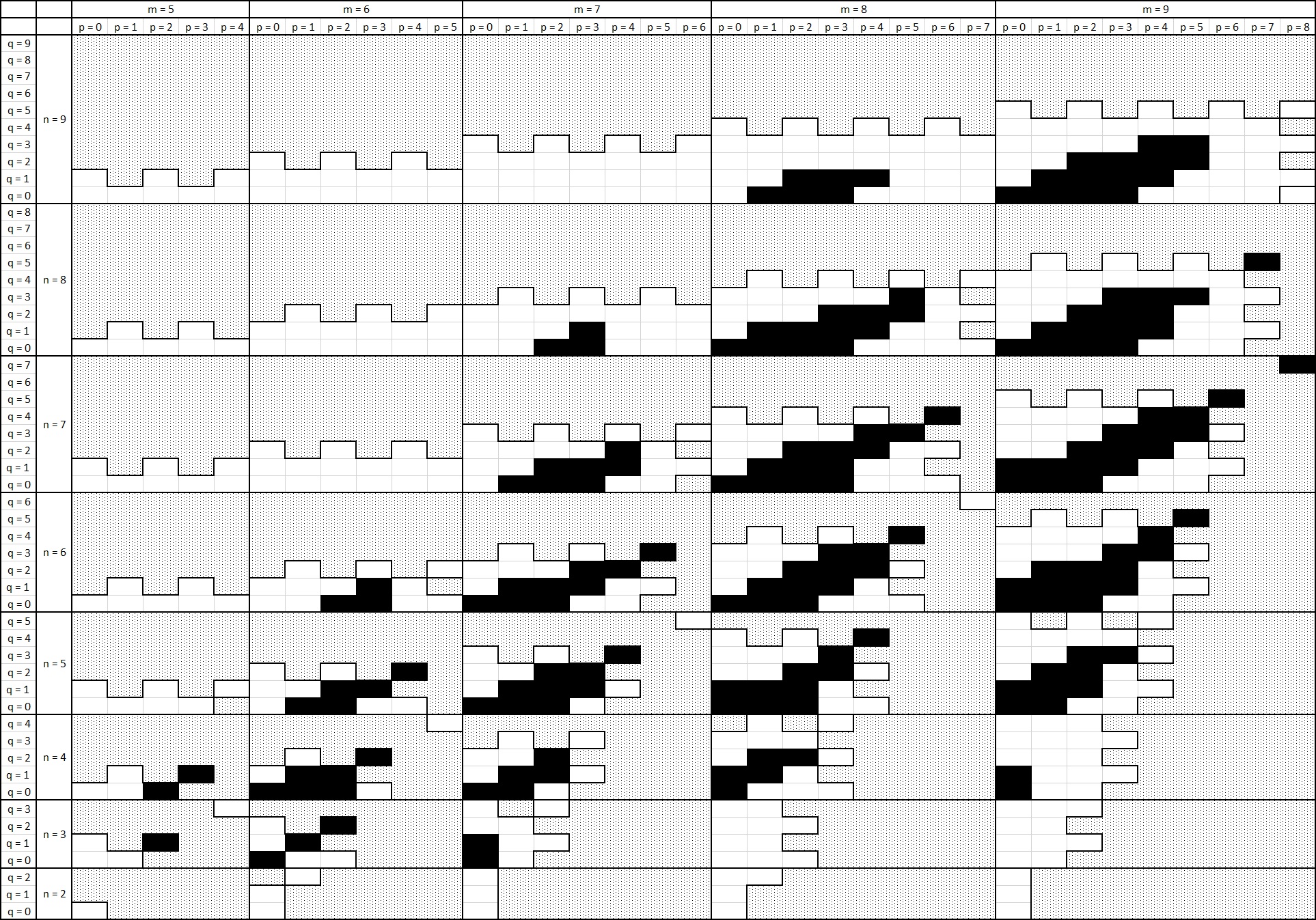}\\
    \label{tableHMon}
    Table \thesubsubsection: Observed monomial ideals of class $\H$.
  \end{center}
\end{minipage}

\hspace{-3.5pc}
\begin{minipage}{35pc}
  \refstepcounter{subsubsection}
  \begin{center}
    \includegraphics[angle=90, width=5.5in]{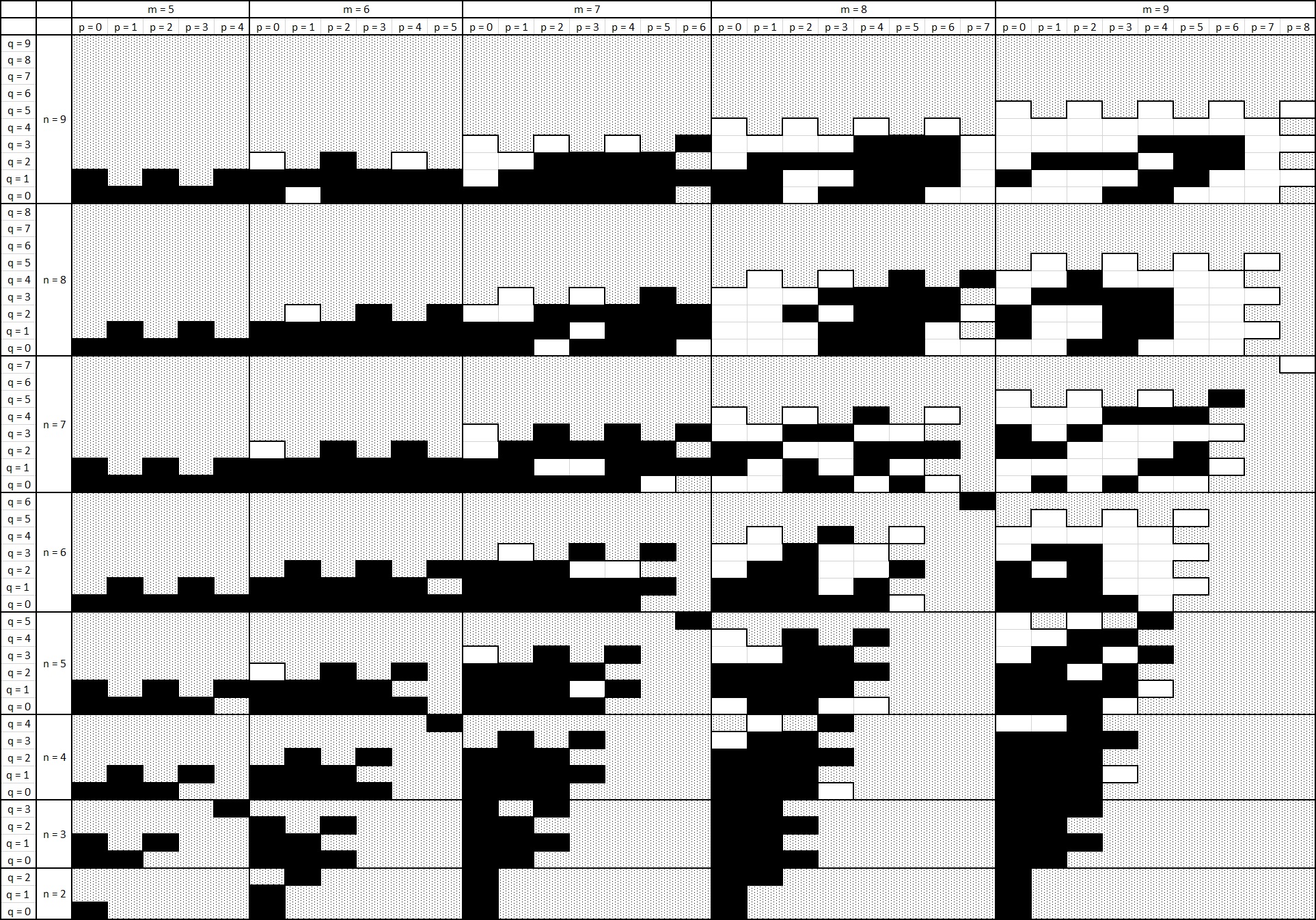}\\
    \label{tableHBin}
    Table \thesubsubsection: Observed binomial ideals of class $\H$.
  \end{center}
\end{minipage}

\subsection{Predominant Classes}
For given values of $m$ and $n$, one may ask if it is possible to
identify a predominant class of ideals, one that is observed more
frequently than others within the experiment. A criterion is needed to
determine when a certain class is observed with such prevalence as to
deserve comment.  In our analysis, this criterion is as follows: in
order to be tagged as predominant, a specific class must have been
observed at least seven times as often as every other class with the
same values of $m$ and $n$. Table \ref{predChar3} has the predominant
classes for $5\leq m\leq 12$ and $2\leq n\leq 10$.  Our criterion does
not identify a predominant class for every pair $(m,n)$. For those
pairs where it fails, all classes that were observed at least
one-seventh as often as the most common class are listed.  In most
cases these lists don't fit in the table and are replaced by labels
explained below.

\smallskip \refstepcounter{subsubsection}
\begin{center}
  {\small \renewcommand{\arraystretch}{1.1}
    \begin{tabular}{|c|c|c|c|c|c|c|c|c|c|c|}
      \hline
      & $m=5$ & $m=6$ & $m=7$ & $m=8$ & $m=9$ & $m=10$ & $m=11$ & $m=12$ \\
      \hline
      $n=10$ & (a) & (j) & (k) & (k) & (u) & (u) & $\H(0,0)$ & $\H(0,0)$ \\
      \hline
      $n=9$ & (b) & (k) & (j) & (k) & (u) & (u) & $\H(0,0)$ & $\H(0,0)$ \\
      \hline
      $n=8$ & (c) & (j) & (k) & (k) & (u) & (u) & $\H(0,0)$ & $\H(0,0)$ \\
      \hline
      $n=7$ & (d) & (l) & (j) & (k) & (u) & $\H(0,0)$ & $\H(0,0)$ & $\H(0,0)$ \\
      \hline
      $n=6$ & (e) & (j) & (k) & (u) & $\H(0,0)$ & $\H(0,0)$ & $\H(0,0)$ & $\H(0,0)$ \\
      \hline
      $n=5$ & (f) & (m) & (q) & $\H(0,0)$ & $\H(0,0)$ & $\H(0,0)$ & $\H(0,0)$ & $\H(0,0)$ \\
      \hline
      $n=4$ & (g) & (n) & (r) & (v) & $\H(0,0)$ & $\H(0,0)$ & $\H(0,0)$ & $\H(0,0)$ \\
      \hline
      $n=3$ & (h) & (o) & (s) & (w) & (v) & $\H(0,0)$ & $\H(0,0)$ & $\H(0,0)$\\
      \hline
      $n=2$ & (i) & (p) & (t) & (x) & (v) & (y) & $\G(3)$ & $\H(0,0)$\\
      \hline 
    \end{tabular}}
  \label{predChar3}
   
  \vspace{.35\baselineskip} Table \thesubsubsection: Predominant
  classes in characteristic $3$.
\end{center}

\smallskip
\begin{center}
  \footnotesize \renewcommand{\arraystretch}{1.1}
  \begin{tabular}{rl|rl}
    (a) & $\H(0,0)$, $\H(1,0)$, $\H(3,0)$, $\H(4,0)$ & (m) & $\H(0,0)$, $\H(1,0)$, $\H(2,0)$, $\H(2,1)$, $\T$ \\
    (b) & $\H(0,0)$, $\H(2,0)$, $\H(3,0)$, $\T$ & (n) & $\H(0,0)$, $\H(1,0)$, $\H(2,1)$ \\
    (c) & $\H(2,0)$, $\T$ & (o) & $\B$, $\H(0,0)$, $\H(2,2)$\\
    (d) & $\H(1,0)$, $\H(2,0)$, $\H(3,0)$, $\T$ & (p) & $\G(3)$, $\H(0,1)$, $\H(1,2)$ \\
    (e) & $\H(0,0)$, $\H(1,0)$, $\H(2,0)$, $\H(3,0)$, & (q) & $\H(0,0)$, $\H(1,0)$, $\H(1,1)$, $\H(2,0)$, \\
        & $\H(4,0)$, $\T$ & & $\H(2,1)$\\
    (f) & $\H(0,0)$, $\H(1,0)$, $\H(2,0)$, $\H(2,1)$, & (r) & $\H(0,0)$, $\H(1,0)$, $\H(1,1)$ \\
        & $\H(3,0)$, $\H(4,1)$, $\T$ & (s) & $\B$, $\H(0,0)$, $\H(0,1)$ \\
    (g) & $\H(1,0)$, $\H(1,1)$, $\H(2,0)$, $\H(3,1)$, $\T$ & (t) & $\B$, $\G(2)$, $\G(4)$, $\H(0,0)$, $\H(0,1)$, $\H(0,2)$ \\
    (h) & $\H(0,1)$, $\H(1,0)$, $\H(2,1)$, $\H(4,3)$ & (u) & $\H(0,0)$, $\H(1,0)$ \\
    (i) & $\B$, $\H(0,0)$ & (v) & $\H(0,0)$, $\H(0,1)$\\
    (j) & $\H(0,0)$, $\H(1,0)$, $\H(2,0)$, $\H(3,0)$ & (w) & $\G(2)$, $\H(0,0)$, $\H(0,1)$ \\
    (k) & $\H(0,0)$, $\H(1,0)$, $\H(2,0)$ & (x) & $\G(3)$, $\G(5)$, $\H(0,0)$ \\
    (l) & $\H(0,0)$, $\H(1,0)$, $\H(2,0)$, $\H(3,0)$, $\T$ & (y) & $\G(2)$, $\G(5)$, $\G(7)$, $\H(0,1)$ \\
  \end{tabular}
\end{center}

The experiment was performed with changing values of $\numTerms$ and
other parameters to seek out ideals of ``rare" classes. To get a
better feeling for the existence of predominant classes, we classified
another $10^6$ randomly generated ideals with $\fieldChar=0$ and all
other parameters set to their default values.  Table \ref{predChar0}
shows the predominant classes in the same fashion as Table
\ref{predChar3}:

\begin{center}
  {\small \renewcommand{\arraystretch}{1.1}
    \refstepcounter{subsubsection}
    \begin{tabular}{|c|c|c|c|c|c|c|c|c|c|c|}
      \hline
      & $m=5$ & $m=6$ & $m=7$ & $m=8$ & $m=9$ & $m=10$ & $m=11$ & $m=12$ \\
      \hline
      $n=10$ & (a) & (j) & (k) & (k) & (l) & (l) & $\H(0,0)$ & $\H(0,0)$ \\
      \hline
      $n=9$ & (b) & (k) & (s) & (k) & (k) & $\H(0,0)$ & (l) & $\H(0,0)$ \\
      \hline
      $n=8$ & (c) & (l) & (q) & (w) & (l) & (l) & $\H(0,0)$ & $\H(0,0)$ \\
      \hline
      $n=7$ & (d) & (m) & $\H(1,0)$ & (q) & $\H(0,0)$ & (l) & $\H(0,0)$ & $\H(0,0)$ \\
      \hline
      $n=6$ & (e) & (n) & $\H(0,0)$ & $\H(1,0)$ & $\H(0,0)$ & $\H(0,0)$ & $\H(0,0)$ & $\H(0,0)$ \\
      \hline
      $n=5$ & (f) & (o) & (t) & $\H(0,0)$ & (\ae) & $\H(0,0)$ & $\H(0,0)$ & $\H(0,0)$ \\
      \hline
      $n=4$ & (g) & (p) & (u) & (x) & $\H(0,0)$ & (\ae) & $\H(0,0)$ & $\H(0,0)$ \\
      \hline
      $n=3$ & (h) & (q) & (v) & (y) & (z) & $\H(0,0)$ & $\H(0,0)$ & (l)\\
      \hline
      $n=2$ & (i) & (r) & $\G(4)$ & (z) & (\o) & (\aa) & none & $\H(0,0)$\\
      \hline
    \end{tabular}
    \label{predChar0}}
    
  \vspace{.35\baselineskip} Table \thesubsubsection: Predominant
  classes in characteristic $0$.
\end{center}

\smallskip

\begin{center}
  \footnotesize \renewcommand{\arraystretch}{1.1}
  \begin{tabular}{rl|rl}
    (a) & $\H(0,0)$, $\H(1,0)$, $\H(4,0)$ & (p) & $\H(0,0)$, $\H(1,0)$, $\H(3,0)$, $\H(5,4)$ \\
    (b) & $\H(3,0)$, $\H(4,0)$, $\T$ & (q) & $\H(0,0)$, $\H(2,0)$ \\
    (c) & $\H(2,0)$, $\H(3,0)$, $\T$ & (r) & $\G(3)$, $\H(0,1)$, $\H(1,2)$\\
    (d) & $\H(1,0)$, $\H(2,0)$, $\H(3,0)$ & (s) & $\H(0,0)$, $\H(1,0)$, $\H(3,0)$, $\T$\\
    (e) & $\H(0,0)$, $\H(1,0)$, $\H(4,0)$, $\T$ & (t) & $\H(0,0)$, $\H(2,0)$,  $\H(2,1)$, $\H(6,5)$\\
    (f) & $\H(0,0)$, $\H(3,0)$, $\H(4,1)$ & (u) & $\H(0,0)$, $\H(1,0)$, $\H(1,1)$ \\
    (g) & $\H(2,0)$, $\H(3,1)$, $\T$ & (v) & $\H(0,0)$, $\H(0,1)$, $\H(2,3)$\\
    (h) & $\H(0,1)$, $\H(1,0)$, $\H(2,1)$, $\H(4,3)$ & (w) & $\H(0,0)$, $\H(1,0)$, $\H(3,0)$ \\
    (i) & $\B$, $\H(0,0)$ & (x) & $\G(2)$, $\H(0,0)$, $\H(2,0)$ \\
    (j) & $\H(0,0)$, $\H(1,0)$, $\H(2,0)$, $\H(3,0)$, $\T$ & (y) & $\G(2)$, $\H(0,0)$, $\H(0,2)$, $\H(1,0)$ \\
    (k) & $\H(0,0)$, $\H(1,0)$, $\H(2,0)$ & (z) & $\G(3)$, $\G(5)$, $\H(0,0)$\\
    (l) & $\H(0,0)$, $\H(1,0)$ &  (\ae) & $\H(0,0)$, $\H(0,1)$, $\H(1,0)$ \\
    (m) & $\H(0,0)$, $\H(2,0)$, $\H(3,0)$, $\H(4,0)$, $\T$ & (\o) & $\H(0,0)$, $\H(0,1)$\\
    (n) & $\H(2,0)$, $\H(3,0)$, $\H(5,2)$, $\T$ & (\aa) & $\G(2)$, $\G(7)$ \\
    (o) & $\H(1,0)$, $\H(2,0)$ & \\
  \end{tabular}
\end{center}


\subsection{Online Repository}\label{Repository}
A text catalogue with examples of ideals from each of the observed
classes can be found on \textit{GitHub}:
\url{https://github.com/ogotchey/codimThreeCode}.  Within the bounds
$5\leq m\leq12$ and $2\leq n\leq 10$ there is a \texttt{*.txt} file
for each class that has been observed. The files are in the
\texttt{data} folder of the repository, and they are named
\texttt{m-n-class-p-q-r.txt} as discussed in Step \ref{step5}. The
lines of each text file are sorted by increasing ``complexity", and
furthermore the files are in \emph{Macaulay2}-readable format.  Also
available in this repository are expanded versions of Tables
\ref{tableBGT}--\ref{tableHBin} and the source code for the algorithm
described in Section \ref{sec1}. For more instructions, see the
$\texttt{README.md}$ file on the repository.

\section{Discussion}

\noindent A main aspect of the realizability question involves bounds
on the parameters $p$, $q$, and $r$ in terms of $m$ and $n$ for ideals
of class $\G$ and $\H$.

\subsection{Ideals of Class $\H$}\label{HBounds}
The current bounds on $p$ and $q$ in terms of $m$ and $n$ for ideals
of class $\H$ were proved in \cite[Theorem 1.1]{CVW20Linkage}. In
Tables \ref{tableH}--\ref{tableHBin} they are represented by the black
border that separates $(p,q)$-cells within each $(m,n)$-box. For
example, the bounds on $p$ and $q$ for $(m,n) = (8,5)$ are seen in the
table below.

\smallskip
\begin{center}
  \refstepcounter{subsubsection}
  \label{tableHBoundEx}
  \includegraphics[height=30mm]{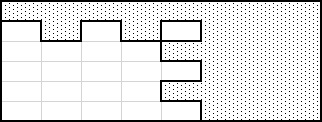}\\
  {Table \thesubsubsection.  Permissible and unrealizable classes for
    $(m,n) = (8,5)$.  }
\end{center}
\smallskip

\noindent
It was conjectured in \cite[Conjecture 7.4]{CVW20Linkage} that for
ideals of class $\H$ with $m\geq 5$ and $n\geq 3$, the bounds
established in \cite[Theorem 1.1]{CVW20Linkage} are optimal. For low
values of $m$ and $n$ we did, indeed, observe ideals of all possible
classes within these bounds; see Table \ref{tableH}. For larger values
of $m$ and $n$ this was not the case. For example, for $m\geq 8$ and
$n\geq 4$, we never observed classes with $p=n-1$ and least possible
$q$, i.e.\ $q=0$ or $q=1$. For $m\geq 8$ and $n\geq 5$, we never
observed classes with $q=m-4$ and least possible $p$, i.e.\ $p=0$ or
$p=1$. For $(m,n) = (8,5)$ these were in fact the only $\H$ classes
not observed, as seen in the table below.

\smallskip
\begin{center}
  \refstepcounter{subsubsection}
  \label{tableRealization}
  \includegraphics[height=30mm]{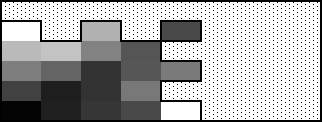}\\
  Table \thesubsubsection. Observed classes for $(m,n) = (8,5)$.
\end{center}
\smallskip

\noindent
While not all permissible classes with $p=n-1$ or $q=m-4$ were
observed in the experiment, it should be pointed out that they are not
unrealizable. Indeed, for $m\ge 5$ and $n\ge 3$, Hardesty proves in
\cite[Theorem 5.2]{Hardesty24} that all these classes can be realized;
the proof is based on linkage theory. Further, within the bounds
$m \le 10$ and $n \le 12$ the experiment established the realizability
of sufficiently many classes to prove, again using linkage, that all
classes within the bounds from \cite[Theorem 1.1]{CVW20Linkage} are
realizable. Thus, the experiment combined with Hardesty's work
provides strong evidence for \cite[Conjecture 7.4]{CVW20Linkage}
regarding ideals of class $\H$.

\subsection{Ideals of Class $\G$} \label{GBounds} The current bounds
on the $r$ parameter for class $\G$, as proved in \cite[Theorem
3.1]{Avramov12}, are represented in Tables
\ref{tableBGT}--\ref{tableBGTBin} by the black borders that separate
$(p,r)$-cells within each $(m,n)$-box. It was conjectured in
\cite[Conjecture 7.4]{CVW20Linkage} that for ideals of class $\G$, if
$n=2$, then $2\leq r\leq m-5$ or $r =m-3$, and if $n\geq 3$, then
$2\leq r\leq m-4$. The observations made in our experiment provide
strong evidence for the conjectured bounds in the $n=2$ case. However,
for $n\geq 3$ the conjectured bound on $r$ seems to be too loose; the
optimal bound on $r$ seems more likely to be a function of both $m$
and $n$---increasing in $m$ and decreasing in~$n$.

\subsection{Classes Realized by Monomial and Binomial Ideals}
As seen in Table \ref{tableHMon}, $\H$ classes realized monomially
tend to have $p$ and $q$ values in close proximity. Classes that can
be realized binomially, see Table \ref{tableHBin}, do not exhibit the
same pattern and, actually, it seems that most $\H$ classes are
realized binomially. One may wonder if every permissible $\H$ class
can be realized by a binomial ideal.

No monomial ideals of class $\G$ were observed in the experiment.
This is further evidence for \cite[Conjecture 6.4]{JLP14} of Painter
as it pertains to class $\G$. (As it pertains to class $\H$, it was
disproved in Faucett's dissertation \cite[Chapter IV]{JAF-phd}.) All
but a few observed $\G$ classes were realized binomially. Similar to
class $\H$, it seems feasible that all permissible $\G$ classes can be
realized binomially.

In the experiment, monomial ideals of class $\B$ were only observed
for $(m,n)$ with $n$ bounded above by $m$ and below by an increasing
function in $m$. On the other hand, binomial ideals of class $\B$ were
observed for all $m \geq 5$ and low values of~$n$.

Monomial ideals of class $\T$ were only observed for $(m,n)$ with $n$
bounded above and below by increasing functions of $m$. There is no
apparent pattern for binomial ideals of class $\T$.

\subsection{Comparing Frequency of the Parameters $p,q,r$}
From Table \ref{tableH} one notices that the most frequently observed
$\H$ classes have values of $p$ and $q$ in close proximity, creating a
gradient that is darkest in the bottom left corner of a given
$(m,n)$-box. That is, the lower values of $p$ and $q$ tend to be
realized more often. This is consistent with Tables \ref{predChar3}
and \ref{predChar0} which show predominant classes.

For class $\G$, low values of $r$ tend to occur more frequently than
high values, see Table \ref{tableBGT}.

\subsection{Comparing Frequency to Linkage} To facilitate this
discussion, we recall that linkage is a symmetric relation on grade
$3$ perfect ideals and say that two classes are \emph{directly linked}
if by one application of linkage one can obtain an ideal of one class
from an ideal of the other. Within the bounds $1\leq n\leq 10$ and
$3\leq m\leq 12$ all permissible $\B$ classes were realized in the
experiment. All but seven permissible $\T$ classes were realized.
However, these seven $\T$ classes are directly linked to observed
classes. In \cite[3.2]{Hardesty24}, Hardesty proves that one can use
linkage to realize ideals of class $\T$ for all $m\geq 5$ and
$n\geq 4$.

Within the bounds of the experiment, $201$ permissible $\G$ classes
were not observed; of these $158$ are directly linked to observed
classes. These $158$ directly linked classes occurred as rarely as
once and as frequently as four million times. Most of these classes
were observed frequently, with all but four of them occurring at least
a hundred times.

A total of $639$ permissible $\H$ classes went unobserved in the
experiment. Of these, $66$ are directly linked to observed
classes. These $66$ directly linked classes occurred as rarely as once
and as frequently as forty thousand times; most of them are rare, with
all but six classes occurring fewer than a hundred times. However,
each of the $639$ permissible $\H$ classes that were not observed is
in the linkage class of an observed class, i.e.\ for each of them
there is a chain of directly linked classes that connects it to an
observed class.

\smallskip

\providecommand{\bysame}{\leavevmode\hbox to3em{\hrulefill}\thinspace}
\providecommand{\MR}{\relax\ifhmode\unskip\space\fi MR }
\providecommand{\MRhref}[2]{%
  \href{http://www.ams.org/mathscinet-getitem?mr=#1}{#2} }
\providecommand{\href}[2]{#2}

\section*{Appendix}
\noindent For given values of $m$ and $n$, it may be useful to have
specific, ``simple" examples of ideals of the various classes. We
provide such examples below; they were first observed in
characteristic $3$ but have been verified to work also in
characteristic $0$.

\definecolor{lgray}{gray}{0.9}
\begin{center}
  \begin{tabular}{|c|c|c|c|c|c|c|}
    \hline 
    \rowcolor{lightgray}
    $m$ & $n$ & class & $p$ & $q$ & $r$ & minimal generators\\ 
    \hline
    $5$ & $2$ & $\B$ & $1$ & $1$ & $2$ & $z^{2},xz,y^{2},xy,x^{2}$ \\
    \rowcolor{lgray}
    $5$ & $2$ & $\H$ & $0$ & $0$ & $0$ & $z^{2},xz,x^{2}y,x^{3}+y^{2}z,y^{4}$ \\
    \hline
    $5$ & $3$ & $\H$ & $0$ & $0$ & $0$ & $yz^{2},x^{2}z,y^{3},x^{2}y+z^{3},x^{3}$ \\
    \rowcolor{lgray}
    $5$ & $3$ & $\H$ & $0$ & $1$ & $1$ & $z^3, yz^2, y^3, xy^2+x^2z, x^3$ \\
    $5$ & $3$ & $\H$ & $1$ & $0$ & $0$ & $xz,yz^{2},y^{3},xy^{2},x^{3}+z^{3}$ \\
    \rowcolor{lgray}
    $5$ & $3$ & $\H$ & $2$ & $1$ & $1$ & $xz,y^{2},z^{3},yz^{2},x^{3}$ \\
    $5$ & $3$ & $\H$ & $4$ & $3$ & $3$ & $z^{2},y^{3},xy^{2},x^{2}y,x^{3}$ \\
    \hline
    \rowcolor{lgray}
    $5$ & $4$ & $\H$ & $0$ & $0$ & $0$ & $z^4, y^2z^2, xy^6-x^6z, x^7, y^8$ \\
    $5$ & $4$ & $\H$ & $1$ & $0$ & $0$ & $xyz+y^{2}z,y^{3},x^{3},z^{4},xz^{3}$ \\
    \rowcolor{lgray}
    $5$ & $4$ & $\H$ & $1$ & $1$ & $1$ & $x^{3}+xyz,z^{4},y^{2}z^{2},x^{2}z^{2},y^{5}$ \\
    $5$ & $4$ & $\H$ & $2$ & $0$ & $0$ & $z^{2},xyz,y^{3},x^{3},x^{2}y^{2}$ \\
    \rowcolor{lgray}
    $5$ & $4$ & $\H$ & $3$ & $1$ & $1$ & $y^{2},z^{3},x^{2}z,x^{3},xyz^{2}$ \\
    $5$ & $4$ & $\T$ & $3$ & $0$ & $0$ & $z^{3},y^{3},x^{3},xyz^{2},xy^{2}z$ \\
    \hline
    \rowcolor{lgray}
    $5$ & $5$ & $\H$ & $0$ & $0$ & $0$ & $xz^{2}+z^{3},y^{2}z,y^{3}-z^{3},x^{2}y,x^{3}-z^{3}$ \\
    $5$ & $5$ & $\H$ & $0$ & $1$ & $1$ & $xy^{2}z,x^{2}yz-z^{4},y^{4}+x^{3}z+y^{2}z^{2},x^{3}y-y^{2}z^{2},x^{4}-z^{4}$ \\
    \rowcolor{lgray}
    $5$ & $5$ & $\H$ & $1$ & $0$ & $0$ & $xz^{2},y^{2}z+z^{3},y^{5}+x^{4}z,x^{6}y,x^{7}+x^{4}y^{3}$ \\
    $5$ & $5$ & $\H$ & $2$ & $0$ & $0$ & $y^{2}z,x^{2}z,y^{3}+z^{3},x^{4}+xyz^{2},x^{3}y^{2}$ \\
    \rowcolor{lgray}
    $5$ & $5$ & $\H$ & $2$ & $1$ & $1$ & $xyz-z^{3},x^{2}y-y^{2}z,z^{5},y^{5}+yz^{4},x^{6}$ \\
    $5$ & $5$ & $\H$ & $3$ & $0$ & $0$ & $xz+z^{2},y^{4}-y^{3}z,xy^{3}-yz^{3},x^{4}-y^{2}z^{2},z^{5}$ \\
    \rowcolor{lgray}
    $5$ & $5$ & $\H$ & $4$ & $1$ & $1$ & $y^{2}+xz,xz^{3}+z^{4},xyz^{2},x^{3}z,x^{4}$ \\
    $5$ & $5$ & $\T$ & $3$ & $0$ & $0$ & $z^{4},xyz^{2},x^{2}y^{2}z,x^{5},y^{6}$ \\
    \rowcolor{lgray}
    \hline
    $6$ & $2$ & $\G$ & $0$ & $1$ & $3$ & $yz,xz,y^{3},xy^{2}+z^{3},x^{2}y,x^{3}$ \\
    $6$ & $2$ & $\H$ & $0$ & $0$ & $0$ & $z^{3},x^{2}z,x^{2}y+y^{2}z,x^{3},y^{4},xy^{3}$ \\
    \rowcolor{lgray}
    $6$ & $2$ & $\H$ & $0$ & $1$ & $1$ & $yz,xz^{2},y^{3},xy^{2}-z^{3},x^{2}y,x^{3}$ \\
    $6$ & $2$ & $\H$ & $1$ & $2$ & $2$ & $xy-z^{2},z^{3},xz^{2},y^{2}z,y^{3},x^{3}$ \\
    \hline
  \end{tabular}
\end{center}

\begin{center}
  \begin{tabular}{|c|c|c|c|c|c|c|}
    \hline
    \rowcolor{lightgray}
    $m$ & $n$ & class & $p$ & $q$ & $r$ & minimal generators\\
    \hline
    \rowcolor{lgray}
    $6$ & $3$ & $\B$ & $1$ & $1$ & $2$ & $xy,x^{2},z^{3},yz^{2},y^{3}z,y^{4}$ \\
    $6$ & $3$ & $\G$ & $0$ & $1$ & $2$ & $z^{3},x^{2}z,y^{3},xy^{2},x^{3}-yz^{2},xyz^{2}$ \\
    \rowcolor{lgray}
    $6$ & $3$ & $\H$ & $0$ & $0$ & $0$ & $yz,xz,xy,z^{3},x^{3},y^{4}$ \\
    $6$ & $3$ & $\H$ & $0$ & $1$ & $1$ & $xz^{2},y^{2}z,x^{2}y,x^{4},xy^{4}-z^{5},y^{6}$ \\
    \rowcolor{lgray}
    $6$ & $3$ & $\H$ & $0$ & $2$ & $2$ & $z^{3},x^{3}z+y^{2}z^{2},y^{4}-x^{2}yz+xy^{2}z,$\\
    \rowcolor{lgray}
        &&&&&& $x^{2}y^{2},x^{3}y-xy^{3},x^{4}+xy^{3}+xy^{2}z$ \\
    $6$ & $3$ & $\H$ & $1$ & $0$ & $0$ & $z^{3},x^{2}z,xy^{2}+yz^{2},x^{3},y^{4}z,y^{5}$ \\
    \rowcolor{lgray}
    $6$ & $3$ & $\H$ & $1$ & $1$ & $1$ & $x^{2},z^{3},yz^{2},xyz,y^{3},xy^{2}$ \\
    $6$ & $3$ & $\H$ & $2$ & $0$ & $0$ & $y^{3},x^{2}y-y^{2}z-xz^{2}-z^{3},x^{3},z^{4},xz^{3},y^{2}z^{2}$ \\
    \rowcolor{lgray}
    $6$ & $3$ & $\H$ & $2$ & $2$ & $2$ & $y^{2},z^{3},xz^{2},x^{2}z,x^{2}y,x^{4}$ \\
    \hline
    $6$ & $4$ & $\B$ & $1$ & $1$ & $2$ & $z^{3},yz^{2},x^{2}y,x^{4},y^{5},xy^{4}z$ \\
    \rowcolor{lgray}
    $6$ & $4$ & $\G$ & $0$ & $1$ & $2$ & $x^{2}z+y^{2}z,x^{3}+xy^{2},yz^{3},y^{3}z,$\\
    \rowcolor{lgray}
        &&&&&& $y^{4}+xyz^{2}-y^{2}z^{2},x^{2}y^{2}-z^{4}$ \\
    $6$ & $4$ & $\H$ & $0$ & $0$ & $0$ & $xz^{2},y^{2}z,x^{2}y,x^{3},z^{4},y^{4}$ \\
    \rowcolor{lgray}
    $6$ & $4$ & $\H$ & $0$ & $1$ & $1$ & $xy^{2}+y^{3},x^{2}y,xz^{3},xyz^{2}+z^{4},x^{3}z-y^{2}z^{2},x^{5}$ \\
    $6$ & $4$ & $\H$ & $1$ & $0$ & $0$ & $xy,z^{3},y^{2}z,x^{2}z,y^{4},x^{4}$ \\
    \rowcolor{lgray}
    $6$ & $4$ & $\H$ & $1$ & $1$ & $1$ & $z^{3},xyz,y^{3},x^{2}y,x^{3},y^{2}z^{2}$ \\
    $6$ & $4$ & $\H$ & $1$ & $2$ & $2$ & $x^{2}z^{3}+xyz^{3},x^{3}z^{2}-x^{2}yz^{2},xy^{3}z+x^{2}yz^{2},$\\
        &&&&&& $x^{2}y^{2}z+x^{2}yz^{2},x^{2}y^{3}+y^{5},x^{5}+y^{5}+x^{2}yz^{2}+z^{5}$ \\
    \rowcolor{lgray}
    $6$ & $4$ & $\H$ & $2$ & $0$ & $0$ & $z^{3},y^{3},xyz^{2},x^{2}yz,x^{2}y^{2},x^{4}$ \\
    $6$ & $4$ & $\H$ & $2$ & $1$ & $1$ & $z^{2},y^{2}z,x^{2}y,x^{3},y^{4},xy^{3}$ \\
    \rowcolor{lgray}
    $6$ & $4$ & $\H$ & $3$ & $0$ & $0$ & 
                                         $y^3z-yz^3, xy^2z+yz^3, y^4-z^4, xz^4+z^5, x^2y^3-z^5, x^5-z^5$\\
    $6$ & $4$ & $\H$ & $3$ & $2$ & $2$ & $y^{2},z^{3},xz^{2},x^{2}z,x^{3}y,x^{4}$ \\
    \rowcolor{lgray}
    $6$ & $4$ & $\H$ & $5$ & $4$ & $4$ & $z^{2},y^{4},xy^{3},x^{2}y^{2},x^{3}y,x^{4}$ \\
    $6$ & $4$ & $\T$ & $3$ & $0$ & $0$ & $xy^{2}z,x^{2}yz,y^{5},x^{5},z^{6},xyz^{5}$ \\
    \hline
    \rowcolor{lgray}
    $6$ & $5$ & $\B$ & $1$ & $1$ & $2$ & $x^{2}z+z^{3},yz^{3}-z^{4},xyz^{2}+z^{4},y^{6}z,x^{4}y^{3},x^{7}+xy^{6}+y^{7}$ \\
    $6$ & $5$ & $\H$ & $0$ & $0$ & $0$ & $z^{3},y^{2}z,x^{2}z+yz^{2},y^{3},xy^{2},x^{3}$ \\
    \rowcolor{lgray}
    $6$ & $5$ & $\H$ & $0$ & $1$ & $1$ & $z^{3},y^{3}z+x^{2}z^{2},xy^{3},x^{2}y^{2}+x^{3}z,x^{3}y+y^{2}z^{2},x^{4}+y^{4}$ \\
    $6$ & $5$ & $\H$ & $1$ & $0$ & $0$ & $y^{3},z^{4},x^{3}z,x^{2}y^{2},xyz^{3},x^{5}$ \\
    \rowcolor{lgray}
    $6$ & $5$ & $\H$ & $1$ & $1$ & $1$ & $y^{3}-xyz,xy^{2},x^{4}z+xyz^{3}-xz^{4},x^{5},xz^{6},z^{8}$ \\
    $6$ & $5$ & $\H$ & $2$ & $0$ & $0$ & $z^{2},xy^{2},x^{3},y^{3}z,x^{2}yz,y^{4}$ \\
    \rowcolor{lgray}
    $6$ & $5$ & $\H$ & $2$ & $1$ & $1$ & $z^{3},yz^{2},x^{3},y^{4},xy^{3},x^{2}y^{2}z$ \\
    $6$ & $5$ & $\H$ & $2$ & $2$ & $2$ & $xy^{2}z+z^{4},xy^{3}-z^{4},x^{3}y-z^{4},xyz^{3}+z^{5},y^{6},x^{6}+x^{3}z^{3}$ \\
    \rowcolor{lgray}
    $6$ & $5$ & $\H$ & $3$ & $0$ & $0$ & $x^{2}y-z^{3},xy^{2}z-z^{4},x^{4}+x^{2}z^{2},yz^{4},y^{2}z^{3},y^{5}+xyz^{3}$ \\
    $6$ & $5$ & $\H$ & $3$ & $1$ & $1$ & $x^{2},z^{4},yz^{3},xy^{2}z,y^{4},y^{3}z^{2}$ \\
    \rowcolor{lgray}
    $6$ & $5$ & $\H$ & $4$ & $0$ & $0$ & $x^{2}y-x^{2}z-xyz,y^{3}z^{2}-x^{2}z^{3}+xyz^{3}+z^{5},$\\
    \rowcolor{lgray}
        &&&&&& $x^{3}z^{2}-xyz^{3}-z^{5},xy^{4}-xyz^{3}+yz^{4}-z^{5},$\\
    \rowcolor{lgray}
        &&&&&& $x^{5}+y^{5}+x^{4}z-x^{2}z^{3}+xyz^{3}-z^{5},yz^{5}+z^{6}$ \\
    $6$ & $5$ & $\H$ & $4$ & $2$ & $2$ & $y^{2},x^{3}z,z^{5},x^{2}z^{3},x^{5},xyz^{4}$ \\
    \rowcolor{lgray}
    $6$ & $5$ & $\T$ & $3$ & $0$ & $0$ & $xy^{2}z,y^{4},x^{2}yz^{2},x^{5},z^{6},xyz^{4}$ \\
    \hline
    $7$ & $2$ & $\B$ & $1$ & $1$ & $2$ & $z^{3},xz^{2},x^{2}z,xy^{2},x^{2}y,x^{3}-xyz,y^{4}$ \\
    \rowcolor{lgray}
    $7$ & $2$ & $\G$ & $0$ & $1$ & $2$ & $z^{3},xz^{2},y^{2}z,xyz,y^{3}+x^{2}z,x^{3}y,x^{5}$ \\
    $7$ & $2$ & $\G$ & $0$ & $1$ & $4$ & $z^{3},y^{2}z,x^{2}z,y^{3},xy^{2},x^{2}y,x^{3}-yz^{2}$ \\
    \rowcolor{lgray}
    $7$ & $2$ & $\H$ & $0$ & $0$ & $0$ & $xyz^{2},x^{3}y,x^{4}z-y^{3}z^{2},y^{5},x^{5},z^{7},xz^{6}$ \\
    $7$ & $2$ & $\H$ & $0$ & $1$ & $1$ & $z^{3},y^{2}z,xy^{2},x^{2}y,x^{2}z^{2},y^{4}+x^{3}z,x^{4}$ \\
    \rowcolor{lgray}
    $7$ & $2$ & $\H$ & $0$ & $2$ & $2$ & $x^{2},xy^{2},y^{2}z^{3},yz^{5},y^{5}z-xz^{5},y^{6},z^{7}$ \\
    \hline
  \end{tabular}
\end{center}

\begin{center}
  \begin{tabular}{|c|c|c|c|c|c|c|}
    \hline 
    \rowcolor{lightgray}
    $m$ & $n$ & class & $p$ & $q$ & $r$ & minimal generators\\
    \hline
    $7$ & $3$ & $\B$ & $1$ & $1$ & $2$ & $z^{3},xz^{2},y^{3},xy^{2},x^{3}z,x^{3}y,x^{4}$ \\
    \rowcolor{lgray}
    $7$ & $3$ & $\G$ & $0$ & $1$ & $2$ & $xz^{2},xy^{2}-z^{3},x^{2}y,x^{3},z^{4},y^{3}z,y^{5}$ \\
    $7$ & $3$ & $\G$ & $0$ & $1$ & $3$ & $yz,xy,z^{4},xz^{3},x^{4}z,y^{5}+x^{3}z^{2},x^{5}$ \\
    \rowcolor{lgray}
    $7$ & $3$ & $\H$ & $0$ & $0$ & $0$ & $z^{3},xz^{2},y^{2}z,xyz,y^{3},x^{2}y,x^{3}$ \\
    $7$ & $3$ & $\H$ & $0$ & $1$ & $1$ & $z^{2},y^{2}z,x^{2}z,y^{3},xy^{2},x^{2}y,x^{4}$ \\
    \rowcolor{lgray}
    $7$ & $3$ & $\H$ & $0$ & $2$ & $2$ & $z^{4},y^{4},xy^{3},x^{2}y^{2},x^{4},x^{2}z^{3},x^{3}yz^{2}+xy^{2}z^{3}-y^{3}z^{3}$ \\
    $7$ & $3$ & $\H$ & $0$ & $3$ & $3$ & $z^{3},xy^{4}-x^{3}yz+x^{2}y^{2}z-xy^{3}z-x^{2}yz^{2},$\\
        &&&&&& $x^{2}y^{3}-x^{4}z-y^{4}z-x^{3}z^{2},y^{5}-x^{4}z-y^{4}z,x^{4}y,$\\
        &&&&&& $x^{3}y^{2}+x^{4}z-x^{3}z^{2}-xy^{2}z^{2},x^{5}-x^{4}z-x^{2}y^{2}z$ \\
    \rowcolor{lgray}
    $7$ & $3$ & $\H$ & $1$ & $0$ & $0$ & $z^{4},y^{4}z,x^{4}z,y^{5}-x^{3}z^{2},x^{4}y,x^{2}y^{4},x^{6}-y^{3}z^{3}$ \\
    $7$ & $3$ & $\H$ & $1$ & $1$ & $1$ & $xy^{2},z^{4},x^{2}z^{2},x^{4}z,x^{4}y,y^{7},x^{9}-y^{6}z^{3}$ \\
    \rowcolor{lgray}
    $7$ & $3$ & $\H$ & $1$ & $2$ & $2$ & $z^{4},xz^{3},x^{2}z^{2},x^{3}z,y^{4},x^{4}+yz^{3},x^{2}y^{3}$ \\
    $7$ & $3$ & $\H$ & $2$ & $1$ & $1$ & $y^{2}z,x^{3},z^{4},yz^{3},x^{2}z^{2}-xyz^{2},x^{2}yz+xz^{3},y^{4}$ \\
    \rowcolor{lgray}
    $7$ & $3$ & $\H$ & $2$ & $3$ & $3$ & $yz^{3},xz^{3},y^{2}z^{2},x^{4},y^{4}z,y^{5},x^{3}y^{3}z-z^{7}$ \\
    \hline
    $7$ & $4$ & $\B$ & $1$ & $1$ & $2$ & $yz^{2},y^{2}z,x^{2}z,y^{3},x^{3},z^{4},xz^{3}$ \\
    \rowcolor{lgray}
    $7$ & $4$ & $\G$ & $0$ & $1$ & $2$ & $x^{2}z,x^{2}y,y^{3}z,y^{4},xy^{3}+z^{4},x^{4}+yz^{3},xy^{2}z^{2}$ \\
    $7$ & $4$ & $\G$ & $0$ & $1$ & $3$ & $xz^{3},xyz^{2},xy^{2}z-y^{3}z-z^{4},y^{4}-x^{2}yz,$\\
        &&&&&& $x^{2}y^{2}-xy^{3}+y^{3}z+z^{4},x^{3}y,x^{4}-y^{3}z+x^{2}z^{2}$ \\
    \rowcolor{lgray}
    $7$ & $4$ & $\H$ & $0$ & $0$ & $0$ & $xz,xy,yz^{2},y^{2}z,x^{3},z^{4},y^{4}$ \\
    $7$ & $4$ & $\H$ & $0$ & $1$ & $1$ & $z^{3},y^{3}+xz^{2},xy^{2},x^{2}y,y^{2}z^{2},x^{3}z,x^{4}$ \\
    \rowcolor{lgray}
    $7$ & $4$ & $\H$ & $1$ & $0$ & $0$ & $y^{3},x^{3},z^{4},xz^{3},y^{2}z^{2},xyz^{2},x^{2}yz$ \\
    $7$ & $4$ & $\H$ & $1$ & $1$ & $1$ & $y^{2},z^{3},xz^{2},xyz,x^{2}z,x^{2}y,x^{3}$ \\
    \rowcolor{lgray}
    $7$ & $4$ & $\H$ & $1$ & $2$ & $2$ & $y^{3},z^{4},xz^{3},xyz^{2},x^{3}z-x^{2}yz-y^{2}z^{2},x^{4}-x^{2}y^{2},x^{2}y^{2}z$ \\
    $7$ & $4$ & $\H$ & $2$ & $1$ & $1$ & $z^{3},y^{3},xy^{2},x^{2}z^{2},x^{3}z,x^{3}y,x^{5}$ \\
    \rowcolor{lgray}
    $7$ & $4$ & $\H$ & $2$ & $2$ & $2$ & $x^{2},xy^{2},z^{4},yz^{3},y^{2}z^{2},y^{4}z,y^{5}$ \\
    $7$ & $4$ & $\H$ & $3$ & $3$ & $3$ & $y^{2},yz^{3},xz^{3},x^{2}z^{2},x^{3}z,x^{4},z^{6}$ \\
    \hline
    \rowcolor{lgray}
    $7$ & $5$ & $\B$ & $1$ & $1$ & $2$ & $y^{2}z,y^{3},x^{3}z,x^{4},z^{5},xz^{4},x^{2}yz^{3}$ \\
    $7$ & $5$ & $\H$ & $0$ & $0$ & $0$ & $z^{3},yz^{2},x^{2}z,y^{3},xy^{2},x^{3}y,x^{4}$ \\
    \rowcolor{lgray}
    $7$ & $5$ & $\H$ & $0$ & $1$ & $1$ & $x^{3},z^{4},yz^{3},y^{2}z^{2},y^{4},xy^{3}-x^{2}yz,x^{2}z^{3}$ \\
    $7$ & $5$ & $\H$ & $1$ & $0$ & $0$ & $z^{2},x^{2}z,xy^{2},y^{3}z,y^{4},x^{3}y,x^{4}$ \\
    \rowcolor{lgray}
    $7$ & $5$ & $\H$ & $1$ & $1$ & $1$ & $z^{3},xz^{2},y^{3},xy^{2}z,x^{3}y,x^{4}z,x^{5}$ \\
    $7$ & $5$ & $\H$ & $2$ & $0$ & $0$ & $z^{3},y^{3},xyz^{2},xy^{2}z,x^{2}yz,x^{4},x^{3}z^{2}$ \\
    \rowcolor{lgray}
    $7$ & $5$ & $\H$ & $2$ & $1$ & $1$ & $z^{2},xy^{2},x^{2}y,y^{3}z,x^{3}z,y^{4},x^{4}$ \\
    $7$ & $5$ & $\H$ & $2$ & $2$ & $2$ & $z^{3},x^{2}y^{2},x^{2}yz^{2},y^{5},xy^{4},x^{5}z,x^{6}$ \\
    \rowcolor{lgray}
    $7$ & $5$ & $\H$ & $3$ & $2$ & $2$ & $z^{3},yz^{2},y^{2}z,x^{4},x^{2}y^{3},y^{6},xy^{5}$ \\
    $7$ & $5$ & $\H$ & $4$ & $3$ & $3$ & $x^{2},z^{4},yz^{3},y^{2}z^{2},y^{3}z,y^{5},xy^{4}$ \\
    \hline
    \rowcolor{lgray}
    $8$ & $2$ & $\G$ & $0$ & $1$ & $2$ & $xyz,x^{2}z,xy^{2},y^{2}z^{2}-xz^{3},x^{4}+y^{3}z,z^{5},yz^{4},y^{5}$ \\
    $8$ & $2$ & $\G$ & $0$ & $1$ & $3$ & $x^{2}z,x^{2}y,y^{3}z,xy^{3},yz^{4},x^{5}-y^{2}z^{3},z^{7},y^{7}+xz^{6}$ \\
    \rowcolor{lgray}
    $8$ & $2$ & $\G$ & $0$ & $1$ & $5$ & $x^{2}z,x^{2}y+xz^{2},x^{3},z^{4},yz^{3},y^{3}z,y^{4},xy^{3}+y^{2}z^{2}$ \\
    $8$ & $2$ & $\H$ & $0$ & $0$ & $0$ & $yz^{2},xz^{2},y^{2}z,x^{2}z,y^{3}+z^{3},xy^{2},x^{2}y,x^{3}-z^{3}$ \\
    \rowcolor{lgray}
    $8$ & $2$ & $\H$ & $0$ & $1$ & $1$ & $xy^{2},z^{4},xz^{3},x^{2}z^{2}-yz^{3},x^{3}z,y^{4},x^{3}y,x^{5}+y^{3}z^{2}$ \\
    $8$ & $2$ & $\H$ & $0$ & $2$ & $2$ & $y^{3},yz^{3},xz^{3},x^{2}z^{2},x^{3}z,x^{2}y^{2}-z^{4},x^{3}y,x^{5}-z^{5}$ \\
    \rowcolor{lgray}
    $8$ & $2$ & $\H$ & $1$ & $2$ & $2$ & $z^{5},x^{4}z,x^{3}y^{2},x^{4}y,y^{5}z,xy^{5},y^{7}+x^{3}z^{4},x^{8}+y^{4}z^{4}$ \\
    \hline
    $8$ & $3$ & $\B$ & $1$ & $1$ & $2$ & $z^{3},y^{2}z^{2},xy^{2}z,xy^{5},x^{4}y^{2},x^{6},y^{7},x^{3}y^{4}+y^{6}z$ \\
    \rowcolor{lgray}
    $8$ & $3$ & $\G$ & $0$ & $1$ & $2$ & $yz^{2},y^{2}z,xyz,x^{2}z,y^{3},xy^{2},x^{2}y-z^{3},x^{4}$ \\
    $8$ & $3$ & $\G$ & $0$ & $1$ & $3$ & $y^{2}z^{3},xy^{4},z^{6},y^{4}z^{2},x^{5}z,y^{6},x^{5}y,x^{6}-yz^{5}$ \\
    \rowcolor{lgray}
    $8$ & $3$ & $\G$ & $0$ & $1$ & $4$ & $z^{3},yz^{2},xyz,y^{4}-x^{3}z,xy^{3},x^{2}y^{2},x^{3}y,x^{4}$ \\
    $8$ & $3$ & $\H$ & $0$ & $0$ & $0$ & $z^{3},yz^{2},y^{2}z,x^{2}z,y^{3}+xz^{2},xy^{2},x^{2}y,x^{3}$ \\
    \rowcolor{lgray}
    $8$ & $3$ & $\H$ & $0$ & $1$ & $1$ & $yz^{2},x^{2}z,y^{3}z,xy^{3},x^{2}y^{2},x^{3}y+z^{4},y^{5},x^{5}$ \\
    $8$ & $3$ & $\H$ & $0$ & $2$ & $2$ & $x^{2},z^{4},yz^{3},xz^{3},y^{2}z^{2},y^{3}z-xyz^{2},y^{4},xy^{3}$ \\
    \hline
  \end{tabular}
\end{center}

\begin{center}
  \begin{tabular}{|c|c|c|c|c|c|c|}
    \hline 
    \rowcolor{lightgray}
    $m$ & $n$ & class & $p$ & $q$ & $r$ & minimal generators\\
    \hline
    \rowcolor{lgray}
    $8$ & $3$ & $\H$ & $1$ & $0$ & $0$ & $y^{3}+z^{3},x^{2}y,z^{4},yz^{3},xyz^{2},xy^{2}z+x^{2}z^{2},x^{3}z,x^{4}-y^{2}z^{2}$ \\
    $8$ & $3$ & $\H$ & $1$ & $1$ & $1$ & $x^{2}y,z^{5},x^{2}z^{3}-y^{2}z^{3},y^{3}z^{2},xy^{3}z,xy^{4},y^{7}-x^{6}z,x^{7}$ \\
    \rowcolor{lgray}
    $8$ & $3$ & $\H$ & $1$ & $2$ & $2$ & $z^{4},yz^{3},xz^{3},x^{2}z^{2},y^{4},x^{4}z,x^{4}y,x^{5}+y^{3}z^{2}$ \\
    $8$ & $3$ & $\H$ & $1$ & $3$ & $3$ & $y^{2},yz^{5},x^{5}z,x^{2}z^{5},x^{6}y-x^{4}yz^{2}-xz^{6},$\\
        &&&&&& $x^{7},z^{8},x^{4}z^{4}+x^{3}yz^{4}$ \\
    \rowcolor{lgray}
    $8$ & $3$ & $\H$ & $2$ & $2$ & $2$ & $xyz+z^{3},yz^{3},xy^{3}-y^{2}z^{2},x^{2}y^{2},x^{3}y,x^{4},y^{4}z,y^{5}$ \\
    \hline
    $8$ & $4$ & $\B$ & $1$ & $1$ & $2$ & $z^{3},yz^{2},y^{2}z,xyz,xy^{2},x^{2}y,x^{3},y^{5}$ \\
    \rowcolor{lgray}
    $8$ & $4$ & $\G$ & $0$ & $1$ & $2$ & $y^{2}z,x^{2}z,x^{2}y,x^{3}-yz^{2},z^{4},xz^{3},y^{4},xy^{3}$ \\
    $8$ & $4$ & $\G$ & $0$ & $1$ & $3$ & $yz^{2},y^{2}z,x^{2}z^{2},y^{5},xy^{4},x^{2}y^{3},x^{4}y+z^{5},x^{6}$ \\
    \rowcolor{lgray}
    $8$ & $4$ & $\G$ & $0$ & $1$ & $4$ & $xy^{3}+y^{4}-y^{3}z,x^{2}y^{2}-y^{3}z,x^{2}yz^{2}+x^{2}z^{3},y^{4}z,$\\
    \rowcolor{lgray}
        &&&&&& $x^{5}-xyz^{3}+xz^{4}+z^{5},y^{2}z^{2},x^{4}y-x^{3}yz-x^{2}z^{3}+yz^{4},$\\
    \rowcolor{lgray}
        &&&&&& $x^{4}z-x^{3}yz-x^{3}z^{2}-x^{2}z^{3}-yz^{4}+z^{5}$ \\
    $8$ & $4$ & $\H$ & $0$ & $0$ & $0$ & $z^{3},yz^{2},xz^{2},y^{2}z,xy^{2},x^{2}y,x^{3},y^{4}$ \\
    \rowcolor{lgray}
    $8$ & $4$ & $\H$ & $0$ & $1$ & $1$ & $z^{2},y^{2}z,x^{3}z,x^{2}y^{2},x^{3}y,y^{5},xy^{4},x^{5}$ \\
    $8$ & $4$ & $\H$ & $0$ & $2$ & $2$ & $x^{3},xz^{4}-yz^{4},y^{2}z^{3},y^{3}z^{2}+xyz^{3},xy^{2}z^{2}+x^{2}z^{3},$\\
        &&&&&& $y^{4}z,y^{5}-xy^{3}z-x^{2}yz^{2}+yz^{4},x^{2}y^{3}+z^{5}$ \\
    \rowcolor{lgray}
    $8$ & $4$ & $\H$ & $1$ & $0$ & $0$ & $xz^{2},yz^{3},y^{4},x^{2}y^{2},x^{4}z,x^{4}y,x^{5}+x^{3}yz+y^{3}z^{2},z^{6}$ \\
    $8$ & $4$ & $\H$ & $1$ & $1$ & $1$ & $x^{3},z^{4},yz^{3},xz^{3},x^{2}z^{2},y^{3}z,x^{2}yz,y^{4}$ \\
    \rowcolor{lgray}
    $8$ & $4$ & $\H$ & $1$ & $2$ & $2$ & $y^{3},yz^{3},xz^{3},x^{2}z^{2},x^{3}z,x^{3}y,z^{6},x^{6}$ \\
    $8$ & $4$ & $\H$ & $2$ & $1$ & $1$ & $x^{2}z^{2},x^{4}+y^{2}z^{2},x^{2}y^{2}z+yz^{4},z^{6},y^{4}z^{2},$\\
        &&&&&& $y^{5}z-xz^{5},xy^{4}z,x^{3}y^{3}+x^{2}y^{4}+xy^{5}+y^{6}$ \\
    \rowcolor{lgray}
    $8$ & $4$ & $\H$ & $2$ & $2$ & $2$ & $y^{3},x^{2}z^{2},x^{2}y^{2},z^{6},xz^{5},x^{6}z,x^{6}y,x^{8}$ \\
    $8$ & $4$ & $\H$ & $3$ & $4$ & $4$ & $z^{3},x^{3}y^{2},y^{6}z,xy^{6},x^{2}y^{5},x^{6}y,x^{9},y^{10}+x^{8}z^{2}$ \\
    \hline
    \rowcolor{lgray}
    $8$ & $5$ & $\B$ & $1$ & $1$ & $2$ & $z^{3},yz^{2},x^{2}yz,x^{3}y,x^{4},xy^{4},y^{5}z,y^{6}$ \\
    $8$ & $5$ & $\G$ & $0$ & $1$ & $2$ & $xz^{3},xy^{3}+z^{4},x^{2}y^{2},x^{3}y+x^{2}z^{2}+xyz^{2},$\\
        &&&&&& $x^{4},y^{3}z^{2},xy^{2}z^{2},y^{6}$ \\
    \rowcolor{lgray}
    $8$ & $5$ & $\H$ & $0$ & $0$ & $0$ & $xy,yz^{2},xz^{2},y^{2}z,x^{2}z,y^{3},x^{3},z^{4}$ \\
    $8$ & $5$ & $\H$ & $0$ & $1$ & $1$ & $x^{3},y^{2}z^{2},xyz^{2},y^{3}z,xy^{3},z^{5},x^{2}z^{4},y^{6}$ \\
    \rowcolor{lgray}
    $8$ & $5$ & $\H$ & $1$ & $0$ & $0$ & $z^{3},xy^{2},x^{2}z^{2},y^{3}z,x^{3}z,y^{4},x^{3}y,x^{4}$ \\
    $8$ & $5$ & $\H$ & $1$ & $1$ & $1$ & $z^{3},xyz,y^{3},xy^{2},x^{2}z^{2},x^{4}z,x^{4}y,x^{5}$ \\
    \rowcolor{lgray}
    $8$ & $5$ & $\H$ & $1$ & $2$ & $2$ & $z^{3},y^{4},xy^{3},x^{2}y^{2},x^{2}yz^{2},$\\
    \rowcolor{lgray}
        &&&&&& $x^{3}y-x^{3}z+x^{2}yz+y^{3}z,x^{4}z,x^{5}$ \\
    $8$ & $5$ & $\H$ & $2$ & $0$ & $0$ & $y^{3},x^{3}z,x^{4}+z^{4},yz^{4},y^{2}z^{3},$\\
        &&&&&& $x^{2}z^{3}+xz^{4},x^{2}yz^{2},x^{2}y^{2}z-xyz^{3}$ \\
    \rowcolor{lgray}
    $8$ & $5$ & $\H$ & $2$ & $1$ & $1$ & $z^{3},yz^{2},x^{3},xy^{2}z,x^{2}y^{2},y^{4}z,y^{5},xy^{4}$ \\
    $8$ & $5$ & $\H$ & $2$ & $2$ & $2$ & $z^{2},y^{3}z,xy^{3},x^{2}y^{2},x^{3}y,x^{5}z,y^{6},x^{6}$ \\
    \rowcolor{lgray}
    $8$ & $5$ & $\H$ & $3$ & $2$ & $2$ & $z^{3},y^{4},xy^{3},x^{2}y^{2},x^{3}z^{2},x^{3}yz,x^{5}y,x^{6}$ \\
    $8$ & $5$ & $\H$ & $3$ & $3$ & $3$ & $y^{3},x^{3}y,z^{5},xz^{4},x^{2}z^{3},x^{3}z^{2},x^{5}z,x^{6}$ \\
    \rowcolor{lgray}
    $8$ & $5$ & $\H$ & $4$ & $4$ & $4$ & $z^{2},x^{2}y^{3},y^{5}z,y^{6},xy^{5},x^{4}y^{2},x^{5}y,x^{6}$ \\
    \hline
  \end{tabular}
\end{center}

\end{document}